\documentclass[a4paper]{cedram-smai-jcm}
\pdfoutput=1
\newif\ifarxiv
\arxivtrue

\ifarxiv
\usepackage{microtype}
\fi

\usepackage[percent]{overpic}
\usepackage[numbers]{natbib}
\usepackage{amsmath}
\usepackage{amscd}
\usepackage{amssymb}
\usepackage{graphicx}
\usepackage{url}
\usepackage{color}
\usepackage{listings}
\usepackage{booktabs}
\usepackage{multirow}
\usepackage{textcomp}
\usepackage{mathtools}
\usepackage{subcaption}
\usepackage{hyperref}
\usepackage[utf8]{inputenc}
\usepackage{tikz}
\usepackage{expl3}
\usetikzlibrary{trees,calc, positioning,decorations.markings,arrows,backgrounds}
\usepackage{pgfplots}
\usetikzlibrary{plotmarks}
\usepackage{pgfplotstable}
\usepackage{subcaption}
\usepgfplotslibrary{groupplots}
\usepackage{environ}
\usepackage{algorithm}
\usepackage{enumerate}
\usepackage[noend]{algpseudocode}
\usepackage{color}
\usepackage{xspace}
\usepackage{svg}
\usepackage{standalone}
\usepackage{xargs}

\makeatletter
\newcommand{\subalign}[1]{%
  \vcenter{%
    \Let@ \restore@math@cr \default@tag
    \baselineskip\fontdimen10 \scriptfont\tw@
    \advance\baselineskip\fontdimen12 \scriptfont\tw@
    \lineskip\thr@@\fontdimen8 \scriptfont\thr@@
    \lineskiplimit\lineskip
    \ialign{\hfil$\m@th\scriptstyle##$&$\m@th\scriptstyle{}##$\crcr
      #1\crcr
    }%
  }
}
\makeatother

\let\div\relax
\DeclareMathOperator{\div}{div}

\DeclareMathOperator{\grad}{grad}
\DeclareMathOperator{\curl}{curl}

\newcommand{\tbnorm}[1]{\vert\kern-0.1em\vert\kern-0.1em\vert\, #1 \,\vert\kern-0.1em\vert\kern-0.1em\vert}

\newcommand{\qbnorm}[1]{\vert\kern-0.15em\vert\kern-0.15em\vert\kern-0.15em\vert\, #1 \,\vert\kern-0.15em\vert\kern-0.15em\vert\kern-0.15em\vert}
\newcommandx{\infsup}[4][1=u,2=v,3=V,4=V]{\underset{\subalign{&{#1} \in {#3} \\ &{#1} \neq 0}}{\mathrm{inf\mathrlap{\phantom{p}}}} \ \underset{\subalign{&{#2} \in {#4} \\ &{#2} \neq 0}}{\mathrm{sup}}\ }

\definecolor{blue}{rgb}{0.2980392156862745, 0.4470588235294118, 0.6901960784313725}
\definecolor{green}{rgb}{0.3333333333333333, 0.6588235294117647, 0.40784313725490196}
\definecolor{red}{rgb}{0.7686274509803922, 0.3058823529411765, 0.3215686274509804}

\renewcommand{\Re}{\ensuremath{\mathrm{Re}}\xspace}
\newcommand{\honev}{\ensuremath{{H}^1(\Omega; \mathbb{R}^d)}\xspace}
\newcommand{\ltwov}{\ensuremath{{L}^2(\Omega; \mathbb{R}^d)}\xspace}
\newcommand{\ltwo}{\ensuremath{{L}^2(\Omega)}\xspace}

\newcommand{\eps}[1]{\ensuremath{\varepsilon(#1)}}

\newcommand{\dx}{\ \mathrm{d}x}

\newcommand{\Pone}{\ensuremath{\mathbb{P}_1}\xspace}
\newcommand{\Ptwo}{\ensuremath{\mathbb{P}_2}\xspace}
\newcommand{\Pk}{\ensuremath{\mathbb{P}_k}\xspace}
\newcommand{\Pkminus}{\ensuremath{\mathbb{P}_{k-1}}\xspace}
\newcommand{\Pkminusdisc}{\ensuremath{\mathbb{P}_{k-1}^\text{disc}}\xspace}
\newcommand{\Pthree}{\ensuremath{\mathbb{P}_3}\xspace}
\newcommand{\PtwoPzero}{\ensuremath{[\mathbb{P}_2]^2\mathrm{-}\mathbb{P}_0}\xspace}

\newcommand{\PtwoPone}{\ensuremath{[\mathbb{P}_2]^2\mathrm{-}\mathbb{P}_1}\xspace}
\newcommand{\PtwoPonedisc}{\ensuremath{[\mathbb{P}_2]^2\mathrm{-}\mathbb{P}_1^\text{disc}}\xspace}
\newcommand{\PthreePtwodisc}{\ensuremath{[\mathbb{P}_3]^2\mathrm{-}\mathbb{P}_2^\text{disc}}\xspace}

\newcommand{\PfivePfourdiscthreed}{\ensuremath{[\mathbb{P}_5]^3\mathrm{-}\mathbb{P}_4^\text{disc}}\xspace}

\newcommand{\advect}[2]{\ensuremath{(#2 \cdot \nabla) #1}}
\newcommand{\mesh}{\ensuremath{\mathcal{M}}\xspace}
\DeclareMathOperator{\MacroStar}{macrostar}
\DeclareMathOperator{\Star}{star}
\newcommand{\dgjmpls}{{\big[\!\!\!\>\big[}}%
\newcommand{\dgjmprs}{{\big]\!\!\!\>\big]}}%
\newcommand{\revnew}[1]{{#1}}%
\usepackage{bm}
\usepackage{xparse}
\usepackage{relsize}
\newcommand{\re}{\mathbb{R}}

\newcommand{\ub}{\mathbf{u}}
\newcommand{\vb}{\mathbf{v}}

\let\temp\phi
\let\phi\varphi
\let\varphi\temp

\newcommand{\vertiii}[1]{{\left\vert\kern-0.25ex\left\vert\kern-0.25ex\left\vert #1 
    \right\vert\kern-0.25ex\right\vert\kern-0.25ex\right\vert}}

\def\Proj_#1{\ensuremath{\operatorname{Proj}_{#1}\!}}

\newcommand{\dr}[1]{\ensuremath{\operatorname{d}\!{#1}}}

\def\drs_#1{\ensuremath{\operatorname{d}_{#1}\!}}
\def\Drs_#1{\ensuremath{\operatorname{D}_{#1}\!}}
\def\Drs^#1{\ensuremath{\operatorname{D}^{#1}\!}}

\def\Ers_#1{\mathrm{E}_{#1}}

\let\div\relax
\DeclareMathOperator{\div}{div}

\newcommand{\defeq}{\vcentcolon=}

\def\xunderbrace#1_#2{{\underbrace{#1}_{\mathclap{#2}}}}
\def\xoverbrace#1^#2{{\overbrace{#1}^{\mathclap{#2}}}}
\def\xunderarrow#1_#2{{\underset{\overset{\uparrow}{\mathclap{#2}}}{#1}}}
\def\xoverarrow#1^#2{{\overset{\underset{\downarrow}{\mathclap{#2}}}{#1}}}

\DeclareDocumentCommand{\boxedeq}{m o}{%
    \IfNoValueTF{#2}{%
        \rlap{\boxed{#1}}%
        \phantom{\hskip\fboxrule\hskip\fboxsep #1}%
        }{%
        \rlap{\boxed{#1#2}}%
        \phantom{\hskip\fboxrule\hskip\fboxsep #1}&\phantom{#2}%
    }%
}

\NewDocumentCommand\meq{m+g}{%
\IfNoValueTF{#2}
  {\overset{\mathsmaller{#1}}{=}}
  {\overset{\mathsmaller{#1}}{\underset{\mathsmaller{#2}}=}}%
}
\catcode`@=11

\def\nvphantom{\v@true\h@false\nph@nt}
\def\nhphantom{\v@false\h@true\nph@nt}
\def\nphantom{\v@true\h@true\nph@nt}
\def\nph@nt{\ifmmode\def\next{\mathpalette\nmathph@nt}%
\else\let\next\nmakeph@nt\fi\next}
\def\nmakeph@nt#1{\setbox\z@\hbox{#1}\nfinph@nt}
\def\nmathph@nt#1#2{\setbox\z@\hbox{$\m@th#1{#2}$}\nfinph@nt}
\def\nfinph@nt{\setbox\tw@\null
    \ifv@ \ht\tw@\ht\z@ \dp\tw@\dp\z@\fi
\ifh@ \wd\tw@-\wd\z@\fi \box\tw@}

\def\XXint#1#2#3{{\setbox0=\hbox{$#1{#2#3}{\int}$ }
\vcenter{\hbox{$#2#3$ }}\kern-.6\wd0}}

\newcommand{\pkg}[1]{\textsf{#1}}

\definecolor{darkblue}{rgb}{0.00,0.00,0.55}
\definecolor{black}{rgb}{0.00,0.00,0.00}
\hypersetup{
    linkcolor = darkblue,
    anchorcolor = darkblue,
    citecolor = darkblue,
    filecolor = darkblue,
    urlcolor = darkblue,
    colorlinks  = true,
}

\title[Reynolds-robust preconditioners for the NSE]{A Reynolds-robust preconditioner for the Scott–Vogelius discretization of the stationary incompressible Navier–Stokes equations}

\author[P.~E.~Farrell]{\firstname{Patrick~E.} \lastname{Farrell}}
\address{Mathematical Institute, University of Oxford, Oxford, UK}
\email{patrick.farrell@maths.ox.ac.uk}

\author[L.~Mitchell]{\firstname{Lawrence} \lastname{Mitchell}}
\address{Department of Computer Science, Durham University, Durham, UK}
\email{lawrence.mitchell@durham.ac.uk}

\author[L.~R.~Scott]{\firstname{L.~Ridgway} \lastname{Scott}}
\address{Department of Computer Science, University of Chicago, Chicago, USA}
\email{ridg@uchicago.edu}

\author[F.~Wechsung]{\firstname{Florian} \lastname{Wechsung}}
\address{Courant Institute of Mathematical Sciences, New York University, New York, USA}
\email{wechsung@nyu.edu}
\thanks{This research is supported by the Engineering and Physical Sciences
Research Council [grant numbers EP/R029423/1 and EP/V001493/1], and by the EPSRC Centre For Doctoral Training in
Industrially Focused Mathematical Modelling [grant number EP/L015803/1] in
collaboration with London Computational Solutions. LM also acknowledges
support from the UK Fluids Network [EPSRC grant number EP/N032861/1] for
funding a visit to Oxford. This work used the ARCHER UK National
Supercomputing Service (\url{http://www.archer.ac.uk}). The authors would
like to acknowledge insightful comments and advice from PD Dr Alexander Linke.}

\begin{document}
\begin{abstract}
Augmented Lagrangian preconditioners have successfully yielded Reynolds-robust
preconditioners for the stationary incompressible Navier--Stokes equations, but
only for specific discretizations. The discretizations for which these
preconditioners have been designed possess error estimates which depend on the
Reynolds number, with the discretization error deteriorating as the Reynolds
number is increased. In this paper we present an augmented
Lagrangian preconditioner for the Scott--Vogelius discretization
on barycentrically-refined meshes. This achieves both Reynolds-robust performance
and Reynolds-robust error estimates. A key consideration is the design of a
suitable space decomposition that captures the kernel of the grad-div term added
to control the Schur complement; the same barycentric refinement that guarantees
inf-sup stability also provides a local decomposition of the kernel of the
divergence. The robustness of the scheme is confirmed by numerical experiments
in two and three dimensions.
\end{abstract}
\maketitle

\ifarxiv
\else

\fi
\section{Introduction} \label{sec:introduction}

The stationary Navier--Stokes equations for the flow of a viscous,
isothermal, incompressible, Newtonian fluid on a bounded Lipschitz domain $\Omega \subset
\mathbb{R}^d$, $d \in \{2, 3\}$, are given by: find $(u, p) \in \honev
\times Q$ such that
\begin{subequations} \label{eq:ns}
  \begin{alignat}{2}
    -  \nabla\cdot 2 \nu \eps{u} + \advect{u}{u} + \nabla p &= f \quad && \text{ in } \Omega, \label{eq:momentum} \\
    \nabla \cdot u &= 0 \quad && \text{ in } \Omega, \\
    u &= g \quad && \text{ on } \Gamma_D, \\
    2\nu \eps{u} \cdot n &= pn \quad && \text{ on } \Gamma_N, \label{eq:naturaloutflow}
  \end{alignat}
\end{subequations}
where $\eps{u} = \frac{1}{2}(\nabla u + \nabla u^\top)$, $\nu > 0$ is the kinematic
viscosity, $f \in \ltwov$, $n$ is the outward-facing unit normal to
$\partial\Omega$, $\Gamma_D$ and $\Gamma_N$ are disjoint with $\Gamma_D \cup
\Gamma_N = \partial \Omega$, and $g \in H^{1/2}(\Gamma_D; \mathbb{R}^d)$. If
$|\Gamma_N| > 0$, then a suitable trial space for the pressure is $Q := \ltwo$;
if $|\Gamma_N| = 0$, then the pressure is only defined up to an additive
constant and $Q := L^2_0(\Omega)$ is used instead~\cite{conca1994}. The Reynolds number
${UL}/{\nu}$, where $U$ and $L$ are the characteristic velocity and length scale
of the flow, is a dimensionless number governing the nature of the system
\cite{landau_fluid_2011}.

For high-Reynolds number flows, it is important that the error estimates
should not degrade as the Reynolds number increases. For most
discretizations, the velocity error estimates are not robust, as they
are polluted by the pressure approximation scaled by the inverse
viscosity. For the Stokes equations, the discrete solution $u_h \in V_h$
can be shown to satisfy
\cite[(3.5)]{john2017}
\begin{equation} \label{eq:errorest}
\|\nabla \left( u - u_h\right)\|_{\ltwo} \le 2 \inf_{\tilde{u}_h \in \mathcal{N}_h} \| \nabla \left( u - \tilde{u}_h\right)\|_{\ltwo} + \nu^{-1} \inf_{q_h \in Q_h} \|p - q_h\|_{\ltwo},
\end{equation}
where $V_h \subset \honev$ is the velocity trial space, $Q_h \subset Q$ is the pressure trial space, and
\begin{equation}
\mathcal{N}_h = \{v_h \in V_h : \int_\Omega q_h\,\nabla \cdot v_h  \dx = 0 \ \forall \ q_h \in Q_h\}
\end{equation}
is the space of discretely divergence-free velocity trial functions.

One way to achieve robustness is to choose a pair $V_h \times Q_h$ such
that $\nabla \cdot V_h \subseteq Q_h$, so that $\nabla \cdot u_h = 0$
holds pointwise \cite{john2017}. If this choice is made, the second
term on the right-hand side of the error estimate \eqref{eq:errorest} for the Stokes equations can be removed and the velocity error is then independent of both pressure and viscosity.
The analysis is more complicated for the Navier--Stokes equations, but
progress has recently been made~\cite{ahmed2018}. We also mention that
similar results are available for the time-dependent case
\cite{linke2019}. Consequently, it is highly advantageous to choose
such a discretization, especially for high Reynolds number flows. Such
discretizations are termed pressure robust or Reynolds-robust.

There are several choices of element pairs that satisfy this
requirement.
They include $H(\div)$-$\ltwo$ discretizations,
such as the Raviart--Thomas and Brezzi--Douglas--Marini families
\cite{cockburn2006,wang2007b,hong2015}; hybrid
discontinuous Galerkin schemes \cite{lehrenfeld2016}; or $H^1$ conforming approaches such as the
Scott--Vogelius $\Pk$-$\Pkminusdisc$ pair \cite{scott_conforming_1985,scott1985}.
The former two options allow for arbitrary order
approximations, but are nonconforming:
the discretization of the momentum equation requires penalty terms. In
contrast, the Scott--Vogelius pair is straightforward to implement, but
is only inf-sup stable on certain types of meshes and for certain polynomial degrees.
In this work, we employ the Scott--Vogelius element.
In~\cite{olshanskii_application_2011}, Olshanskii \& Rebholz demonstrate the accuracy of this discretization and investigate the numerical performance of direct sparse solvers.
While these solvers perform very well for problems with sizes of the order of
millions of degrees of freedom, they do not scale and a different strategy is
required to solve larger problems.

In this work, we build on the insights of
\cite{benzi2006,hong2015,farrell2018b, fmsw_robust_2020} to develop a Reynolds-robust
block preconditioner for the Reynolds-robust Scott--Vogelius discretization of the
Navier--Stokes equations.
When building block preconditioners, the usual difficulty is
developing a good approximation for the inverse of the Schur complement.
We employ an augmented Lagrangian term
to control the Schur complement of the system; this simplifies
the approximation of the Schur complement, at the cost of making
the momentum equation significantly more difficult to solve.
We then develop
a specialized geometric multigrid
scheme for the resulting augmented momentum equation.
This approach is based on the work of Sch\"oberl \cite{schoberl1999b}. Its application to the Navier--Stokes equations was pioneered in two dimensions by Benzi \& Olshanskii~\cite{benzi2006},
and has recently been extended to three dimensions~\cite{farrell2018b}.
A similar strategy has proven successful for a $H(\div)$-$\ltwo$ discretization of the Stokes
equations \cite{hong2015}.
While effective at controlling
iteration counts as the Reynolds number is varied, the schemes presented in \cite{benzi2006,farrell2018b} heavily rely on the use of piecewise constant pressure functions and are not effective for Scott--Vogelius discretizations.

Many alternative approaches to solving \eqref{eq:ns} have been considered in
the literature. These include the pressure convection-diffusion (PCD),
least-squares commutator (LSC) and SIMPLE block preconditioners \cite{patankar1980,kay2002,elman2006,elman2014}, monolithic
multigrid approaches \cite{vanka1986,turek1999}, and a modified augmented Lagrangian
strategy which trades off control of the Schur complement against ease of
solving the augmented momentum block \cite{benzi2011}. These preconditioners do not generally
enjoy Reynolds-robustness.

The remainder of this manuscript is structured as follows.
In section \ref{sec:sv-stability} we discuss the discretization and the conditions under which inf-sup stability of the Scott--Vogelius element is known.
In section \ref{sec:augmented_lagrangian} we recall the augmented Lagrangian
preconditioning strategy and the difficulties it introduces for solving the momentum equations.
A Reynolds-robust
multigrid cycle for the augmented Scott--Vogelius momentum operator is introduced
in section \ref{sec:multigrid}. Finally, numerical examples in two and three dimensions
are presented in section \ref{sec:examples}.
\section{Discretization}\label{sec:sv-stability}
We begin by recalling that finite dimensional subspaces $V_h\subset V$ and $Q_h\subset Q$ are said to satisfy the inf-sup condition if there exists a $\gamma>0$ such that
\begin{equation}
    \infsup[q_h][v_h][Q_h][V_h] \frac{(\nabla\cdot v_h, q_h)}{\|v_h\| \|q_h\|} \ge \gamma
\end{equation}
for all mesh sizes $h>0$.
Intuitively, the inf-sup condition encourages large velocity and small pressure spaces, but the opposite is true for the condition $\nabla\cdot V_h \subseteq Q_h$.
For this reason it is difficult to construct discretizations that satisfy both of these properties. Most of the popular inf-sup stable finite element discretizations of the Stokes and Navier--Stokes equations (such as the Taylor--Hood, the MINI or the $\PtwoPzero$ elements) do not satisfy $\nabla \cdot V_h\subseteq Q_h$.

The Scott--Vogelius element is given by choosing continuous piecewise polynomials of degree $k$ for the velocity and discontinuous piecewise polynomials of degree $k-1$ for the pressure.
While this clearly implies that $\nabla \cdot V_h\subseteq Q_h$, inf-sup stability of the Scott--Vogelius element is more delicate, and is a topic of ongoing research.
In two dimensions, Scott \& Vogelius proved~\cite{scott1985} that the element is inf-sup stable for $k\ge 4$ if the mesh does not have nearly singular vertices.
In three dimensions, it was proven more recently in \cite{zhang2011} that the element is stable for $k\ge 6$ on uniform meshes.
The stability on general tetrahedral meshes continues to be an open question \cite{neilan_2020}.

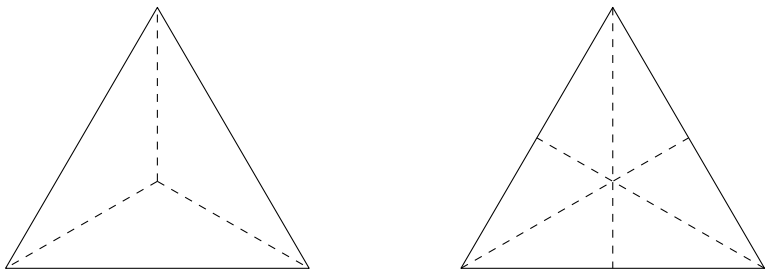
\begin{figure}[htbp]
  \centering
  \ExplSyntaxOn
\newcommand{\fpeval}[1]{\fp_eval:n{#1}}
\ExplSyntaxOff

\begin{tikzpicture}[scale=4]
    \begin{scope}[shift={(0,0)}]
        \coordinate (0) at (0, 0);
        \coordinate (1) at (1, 0);
        \coordinate (2) at (0.5, \fpeval{sqrt(3)/2});
        \coordinate (bc) at (barycentric cs:0=1,1=1,2=1);
        \draw[line join=round] (0) -- (1) -- (2) -- cycle;
        \draw[line join=round, dashed] (bc) -- (0);
        \draw[line join=round, dashed] (bc) -- (1);
        \draw[line join=round, dashed] (bc) -- (2);
    \end{scope}
    \begin{scope}[shift={(1.5,0)}]
        \coordinate (0) at (0, 0);
        \coordinate (1) at (1, 0);
        \coordinate (2) at (0.5, \fpeval{sqrt(3)/2});
        \draw[line join=round] (0) -- (1) -- (2) -- cycle;
        \draw[line join=round, dashed] (0) -- (barycentric cs:1=1,2=1);
        \draw[line join=round, dashed] (1) -- (barycentric cs:0=1,2=1);
        \draw[line join=round, dashed] (2) -- (barycentric cs:0=1,1=1); \coordinate (0) at (0, 0);
    \end{scope}
\end{tikzpicture}
  \caption{Barycentrically refined triangle (also known as Alfeld split) on the left, and Powell--Sabin split on the right.\label{fig:alfeld-powell-sabin}}
\end{figure}

On barycentrically refined meshes, however, the pair is known to be stable for polynomial order $k = d$, see~\cite[Section 4.6]{qin1994} for the 2D case and~\cite{zhang2004} for the 3D case.
If one is willing to consider the more complicated Powell--Sabin split, the order can be reduced further to $k=d-1$ \cite{zhang_p1_2008,zhang_quadratic_2011}.
The two refinement patterns are shown for the two dimensional case in Figure~\ref{fig:alfeld-powell-sabin}.
In this work we will consider the case of $k \ge d$ on barycentrically refined meshes, but the arguments apply \emph{mutatis mutandis} to the Powell--Sabin split.

In the context of the multigrid scheme that we will develop in section~\ref{sec:multigrid}, the requirement for barycentrically refined elements has some implications for our mesh hierarchy.
First, note that repeatedly barycentrically refining a mesh leads to degenerate elements.
Furthermore, it is not known whether regularly refining a mesh, on which an element pair with $k=d$ is stable, always results in a refined mesh for which stability is maintained.
Consequently, we build the multigrid hierarchy in a different way.
Given a domain $\Omega$, we consider a given simplicial mesh $\mesh_h = \{K^h\}$ with $\cup_{K^h\in\mesh_h} K^h = \overline{\Omega}$ and $(K^h_i)^\circ \cap (K^h_j)^\circ = \emptyset$ for $i \neq j$.
The elements $K^h\in\mesh_h$ will be referred to as the \emph{macro cells}.
For each level $h$, we obtain the mesh $\hat\mesh_h$ by barycentric refinement: that is, for each macro cell $K^h\in\mesh_h$ we obtain $d+1$ many cells $\hat K^h_{i}$, $0\le i\le d$ and
\begin{equation}
    \hat\mesh_h = \{\hat K^h_{i}: 0\le i\le d, K^h\in\mesh_h\}.
\end{equation}
The function spaces on $\hat \mesh_h$ are then given by
\begin{align}
    V_h &\defeq \{ \vb \in H^1(\Omega;\re^d): \vb\vert_{\hat K} \in [\mathbb{P}^d(\hat K)]^{d} \ \forall \hat K \in \hat\mesh_h\},\\
    Q_h &\defeq \{ q \in L^2(\Omega): q\vert_{\hat K} \in \mathbb{P}^{d-1}(\hat K) \ \forall \hat K \in \hat\mesh_h\}.
\end{align}

\begin{figure}[htbp]
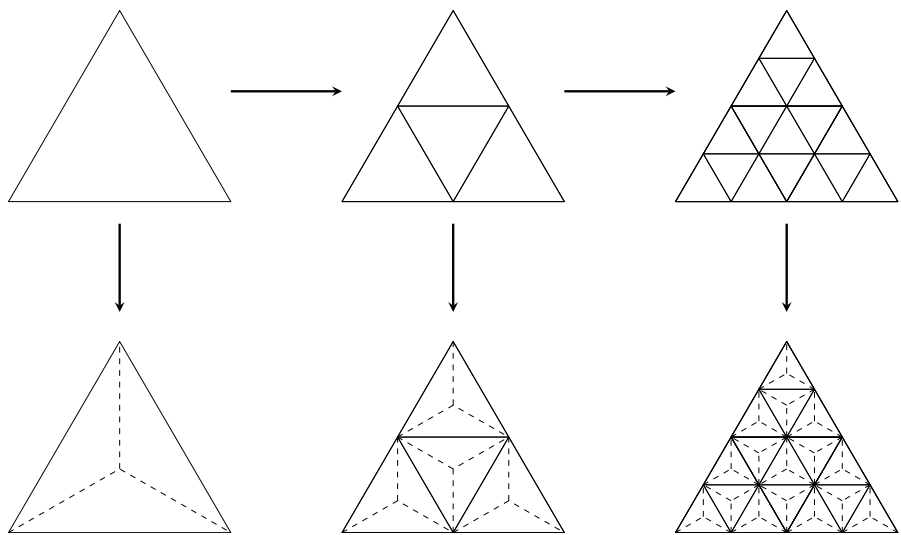

  \centering
  \includestandalone[width=0.7\textwidth]{./figures/bary-hierarchy}
  \caption{A three level barycentrically refined multigrid hierarchy.}\label{fig:bary-hierarchy}
\end{figure}

We construct the hierarchy as follows. We start with an initial coarse triangulation of the domain, given by $\mesh_H$.
We obtain $\mesh_h$, $h=\frac{1}{2}H$, by uniform refinement of the initial mesh.
Both $\mesh_H$ and $\mesh_h$ are then refined barycentrically to obtain $\hat\mesh_H$ and $\hat\mesh_h$.
Note that though $\mesh_H$ and $\mesh_h$ form a nested hierarchy, this is not true for $\hat\mesh_H$ and $\hat\mesh_h$.
This two-level approach canonically extends to many levels; a hierarchy of three levels is shown in Figure~\ref{fig:bary-hierarchy}.

We will see in Section~\ref{ssec:relaxation} that this macro element structure not only guarantees inf-sup stability, but is also crucial in defining a robust relaxation method.


\section{Variational formulation and augmented Lagrangian strategy} \label{sec:augmented_lagrangian}
For boundary data $g \in H^{1/2}(\Gamma_D)$, let
\begin{equation}
V_g = \{v \in H^1(\Omega; \mathbb{R}^d) : \left. v \right|_{\Gamma_D} = g \}.
\end{equation}
The weak form of \eqref{eq:ns} is: find $(u, p) \in V_g \times Q$ such that
\begin{equation}
    \int_\Omega 2\nu\eps{u}: \eps{v} \dx + \int_\Omega \advect{u}{u} \cdot v \dx - \int_\Omega p \nabla \cdot v \dx - \int_\Omega q \nabla \cdot u \dx = \int_\Omega f \cdot v \dx,
\end{equation}
for all $(v, q) \in V_0 \times Q$.

Given finite dimensional subspaces $V_h\subset \honev$ and $Q_h \subset Q$ and after applying Newton's method to the nonlinear equations, at each Newton step we must solve a nonsymmetric linear system of generalized saddle point structure:
\begin{equation}
  \label{eq:ns-saddle-point}
  \begin{pmatrix}
    A & B^\top\\
    B & 0
  \end{pmatrix}
  \begin{pmatrix}
    \delta u\\
    \delta p
  \end{pmatrix}
  =
  \begin{pmatrix}
    b\\
    c
  \end{pmatrix},
\end{equation}
where $A$ is the discrete linearized momentum operator, $B^\top$ the
discrete gradient, $B$ the discrete divergence, and $\delta u$ and
$\delta p$ are the updates to the velocity and pressure solutions
respectively.
This system becomes increasingly difficult to solve as the
Reynolds number is increased.

There are two key ingredients when building block preconditioners for \eqref{eq:ns-saddle-point}:
an effective solver for $A$ and an effective solver for the Schur complement $S = -B A^{-1} B^\top$.  Since $S$ is usually dense, tractable approximations
to $S^{-1}$ must be developed on a PDE-specific basis.
The main issue with PCD, LSC and SIMPLE is that their
choice for the approximate Schur complement becomes a poor
approximation to the true Schur complement as the Reynolds number is increased,
which in turn results in significant growth of the iteration counts.
In \cite{benzi2006}, Benzi \& Olshanskii proposed an augmented Lagrangian strategy that
significantly simplifies the approximation of the Schur complement:
for $\gamma>0$, the linear system  \eqref{eq:ns-saddle-point} is augmented by
adding a term to the top-left block and adjusting the residual accordingly:
\begin{equation}
    \begin{bmatrix}
        A + \gamma B^\top M_p^{-1} B & B^\top \\
        B & 0
    \end{bmatrix}
    \begin{bmatrix}
        \delta u \\ \delta p
    \end{bmatrix}
    =
    \begin{bmatrix}
        b + \gamma B^\top M_p^{-1} c\\ c
    \end{bmatrix},
\end{equation}
where $M_p$ is the mass matrix for the pressure space.
It is immediately clear that this modification does not change the solution of the linear system.
Furthermore, one can show that the Schur complement of the augmented system, $\tilde S$, satisfies
\begin{equation}\label{eqn:schur-complement}
    \tilde S^{-1} = S^{-1} - \gamma M_p^{-1},
\end{equation}
where $S$ is the Schur complement of the original system.
The advantage is clear: as $\gamma\to\infty$, $\tilde S^{-1}$ can be
approximated by a scaled inverse pressure mass matrix, which is easy to solve.

In general, a triple matrix product as it occurs in the augmented Lagrangian term is both expensive to compute and store.
However, it is straightforward to check that adding $\gamma B^\top M_p^{-1} B$ to the linear system corresponds to augmenting the weak form with a term 
\begin{equation}\label{eqn:al-weak-form}
    \gamma \int_\Omega \Pi_{Q_h} (\nabla\cdot u) \Pi_{Q_h} (\nabla\cdot v)\dr x,
\end{equation}
where $\Pi_{Q_h}$ is the projection onto $Q_h$.

The same augmentation, without the
projection onto $Q_h$, is known as grad-div stabilization as it corresponds to the weak form  of
$-\gamma \nabla\nabla\cdot u$.
\revnew{We remark that this form does not lead to any additional non-zeros in the assembled linear system.}
As we focus on discretizations that satisfy $\nabla\cdot V_h\subseteq Q_h$, in this work the projection onto the pressure space is always the identity operator and hence the augmented Lagrangian and grad-div stabilization coincide.

To summarize, the augmented problem in weak form reads: 
find $(u, p) \in \left(V_h \cap V_g\right) \times Q_h$ such that
\begin{equation}
    \begin{aligned}
        &\int_\Omega 2\nu\eps{u}: \eps{v} \dx + \int_\Omega \advect{u}{u} \cdot v \dx + \gamma \int_\Omega \nabla\cdot u\nabla\cdot v\dx\\
        &\qquad\qquad\qquad\qquad    - \int_\Omega p \nabla \cdot v \dx - \int_\Omega q \nabla \cdot u \dx = \int_\Omega f \cdot v \dx,
    \end{aligned}
\end{equation}
for all $(v, q) \in \left(V_h \cap V_0\right) \times Q_h$.
\revnew{
\begin{remark}\label{rem:th-schur}
A discretization for which the augmented Lagrangian and the grad-div stabilization \emph{do not} coincide is the classical $[\Pk]^d\mathrm{-}\Pkminus$ Taylor--Hood element.
While the Taylor--Hood element does not enforce the divergence constraint exactly, it can be shown that when using grad-div stabilization, in the limit of $\gamma\to\infty$ the solution converges to that obtained when using the Scott--Vogelius discretization, so long as the Scott--Vogelius element is stable \cite[Theorem 3.1]{case2011}.

    It turns out that not only are the solutions obtained from these two elements related, but also that (on the particular meshes considered here) the solver developed for the Scott--Vogelius element can be used for the Taylor--Hood element.
    First observe that since the pressure space is continuous, the projection onto $Q_h$ in the augmented Lagrangian term introduces additional coupling and adding \eqref{eqn:al-weak-form} increases the number of non-zeros in the assembled matrix.
However, it was shown in~\cite[Lemma~2]{heister2012} that in the limit of $h\to 0$ grad-div stabilization converges to the discrete augmented Lagrangian approach.
    Hence one can simply use grad-div stabilization and still use the Schur complement approximation in~\eqref{eqn:schur-complement}.
    Clearly, the obtained top-left block is then identical to that obtained from the Scott--Vogelius element and hence the multigrid scheme that we propose in the next section applies in the same way.
    From a solver point of view, a minor disadvantage of this approach is that the pressure mass matrix corresponding to a continuous finite element space is more expensive to invert than that of a discontinuous space.
\end{remark}
}

\section{Solving the top-left block} \label{sec:multigrid}
Although the augmented Lagrangian approach is appealing since it provides
for an excellent Schur complement approximation, it introduces some new
challenges to developing fast solvers for the equations as a whole. In contrast to $A$,
\revnew{the performance of standard multigrid schemes deteriorates quickly for $A_\gamma = A + \gamma
B^{\top}M_p^{-1}B$ as $\gamma$ is increased.} The additional term has a large kernel consisting of
solenoidal vector-fields, implying that the problem becomes nearly singular as
$\gamma$ increases. If a scalable solver for the nearly singular top-left block
can be developed that is robust to $\nu$ and $\gamma$, the overall solver
will be Reynolds-robust.

Multigrid schemes for nearly singular problems have been studied previously by
Sch\"oberl \cite{schoberl1999,schoberl1999b, schoberl_robust_1998} and Lee et
al.~\cite{lee2007,lee2009}, and a similar analysis was carried out for
overlapping Schwarz methods by Ewing and Wang~\cite{ewing1992}. \revnew{While these works only consider the symmetric case, they clearly demonstrate that a necessary condition for} a $\gamma$-robust scheme is a good understanding and
characterization of the kernel of the semidefinite term, which in this case
means understanding the kernel of the divergence operator.
This has implications for both multigrid relaxation and prolongation, as we
now consider.

\subsection{Relaxation} \label{ssec:relaxation}
The core requirement for obtaining a parameter robust relaxation is that the
space decomposition defining the relaxation needs to provide a decomposition of
the kernel of the singular operator.
Classical point relaxation methods such as point-block Jacobi do not satisfy
this property and are ineffective
for this problem; their smoothing strength degrades as $\gamma$ increases.

Many smoothers can be expressed as so-called subspace correction methods
\cite{xu1992}.
We consider a decomposition
\begin{equation}
    V_h = \sum_i V_i,
\end{equation}
where the sum is not necessarily a direct sum. This decomposition naturally
defines an associated additive Schwarz relaxation method.
For each subspace $i$ we then denote the natural inclusion by $I_i:V_i \to V_h$
and we define the restriction $A_i$ of $A$ onto $V_i$ as 
\begin{equation}
    (A_i u_i, v_i) = (A I_i u_i, I_i v_i) \qquad \text{ for all } u_i,\ v_i\in V_i.
\end{equation}
Now denoting
\begin{equation}
    D_h^{-1} = \sum_i I_i A_i^{-1} I_i^\top
\end{equation}
and introducing a damping parameter $\tau>0$, we can express one update of the additive Schwarz method as
\begin{equation}
    u_{k+1}  = u_{k} + \tau D^{-1}(f-Au_k).
\end{equation}
The method is also known as the parallel subspace correction method \cite{xu1992}.
For the choice $V_i = \operatorname{span}\big(\{\phi_i\}\big)$, we recover the classical Jacobi iteration.
Usually each subspace is described by an index set $J_i$ and $V_i = \operatorname{span}\big(\{\phi_j: j\in J_i\}\big)$.
In that case we also speak of a block Jacobi method.

Both Sch\"oberl \revnew{(for additive relaxation)} and \revnew{later} Lee et al.~\revnew{(for multiplicative relaxation)} recognize that the key condition for a $\gamma$-robust relaxation is that the subspaces $V_i$ need to satisfy a \emph{kernel decomposition property}:
\begin{equation}\label{eqn:kernel-decomposition}
    \mathcal{N}_h = \sum_i (V_i\cap \mathcal{N}_h),
\end{equation}
where we recall that $\mathcal{N}_h$ is the space of discretely divergence-free vector-fields.
In essence, we require that any kernel function can be written as the sum of kernel functions
drawn from the subspaces.

\begin{figure}[htbp]
    \centering
    \begin{tikzpicture}[scale=4,
    bcdof/.style={postaction={decorate,decoration={markings,
                mark=at position 0.5 with {\draw[solid, fill, black!60] circle
        (0.12);}}}},
        dof/.style={postaction={decorate,decoration={markings,
            mark=at position 0.5 with {\draw[solid, fill, black] circle (0.24);}}}
            }]

            \foreach \i in {0, 1, 2, 3, 4} {
                \coordinate (0-\i) at ($(\i*1.5, 0)$);
                \coordinate (1-\i) at ($(0-\i) + (60:1.5)$);
                \coordinate (2-\i) at ($(1-\i) + (120:1.5)$);
            }
            \foreach \j in {0, 2} {
                \foreach \i in {1, 2, 3, 4} {
                    \draw[line width=1mm, solid, black, fill] (\j-\i) circle (0.06);
                }
            }
            \foreach \i in {0, 1, 2, 3, 4} {
                \draw[line width=1mm, solid, black, fill] (1-\i) circle (0.06);
            }

            \foreach \i/\j in {1/2, 2/3, 3/4} {
              \foreach \k in {0, 2} {
                \draw[dof, line join=miter, line width=1.0mm] (\k-\i) -- (\k-\j);
              }
              \draw[dof, line join=miter, line width=1.0mm] (0-\i) -- (1-\i);
              \draw[dof, line join=miter, line width=1.0mm] (1-\i) -- (2-\i);
              \draw[dof, line join=miter, line width=1.0mm] (0-\j) -- (1-\i);
              \draw[dof, line join=miter, line width=1.0mm] (1-\i) -- (2-\j);
            }

           \draw[dof, line join=miter, line width=1.0mm] (1-0) -- (2-1);
           \draw[dof, line join=miter, line width=1.0mm] (1-0) -- (0-1);
           \draw[dof, line join=miter, line width=1.0mm] (1-4) -- (0-4);
           \draw[dof, line join=miter, line width=1.0mm] (1-4) -- (2-4);
            \foreach \i/\j in {0/1, 1/2, 2/3, 3/4} {
              \draw[dof, line join=miter, line width=1.0mm] (1-\i) -- (1-\j);
            }

  \begin{scope}[on background layer]
            \draw[blue!50, fill=blue!50] 
                ([shift={(0:10pt)}]1-0) --
                ([shift={(60:10pt)}]0-1) --
                ([shift={(120:10pt)}]0-2) --
                ([shift={(180:10pt)}]1-2) --
                ([shift={(240:10pt)}]2-2) --
                ([shift={(300:10pt)}]2-1) --
                cycle;
            \draw[blue!50, fill=blue!50]
                ([shift={(0:10pt)}]1-2) --
                ([shift={(60:10pt)}]0-3) --
                ([shift={(120:10pt)}]0-4) --
                ([shift={(180:10pt)}]1-4) --
                ([shift={(240:10pt)}]2-4) --
                ([shift={(300:10pt)}]2-3) --
                cycle;
  \end{scope}





\end{tikzpicture}
    \caption{The star patch satisfies the kernel decomposition property \eqref{eqn:kernel-decomposition} for the $\PtwoPzero$ element.}
    \label{fig:star}
\end{figure}
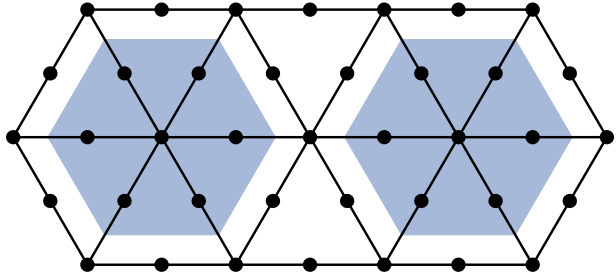

When considering the $\PtwoPzero$ element on a 2D triangulation $\mesh_h$, that is
\begin{equation}
    \begin{aligned}
        V_h &= \{ \vb \in H^1(\Omega;\re^d): \vb\vert_{K} \in [\mathbb{P}^2(K)]^{d} \ \forall K \in \mesh_h\},\\
        Q_h &= \{ q \in L^2(\Omega): q\vert_{K} \equiv \text{const} \ \forall K \in \mesh_h\},
    \end{aligned}
\end{equation}
Sch\"oberl proved that the space decomposition defined by
\begin{equation}
    V_i \coloneqq \{ \ub_h \in V_{h} : \mathrm{supp}(\ub_h) \subset \Star(v_i) \},
\end{equation}
where
\begin{equation}
    \Star(v_i) \coloneqq \bigcup_{{K \in \mesh_h \,: \,v_i \in K}} K,
\end{equation}
satisfies the kernel decomposition property~\eqref{eqn:kernel-decomposition}.
The subspaces are illustrated in Figure~\ref{fig:star}.
This decomposition was subsequently used in \cite{benzi2006} for the
Navier--Stokes equations in two dimensions and in \cite{farrell2018b} in three
dimensions. It was also used for $H(\div)$-$\ltwo$ discretizations of the Stokes
and linear elasticity equations in \cite{hong2015}, and provides a robust
relaxation method for the $H(\div)$ and $H(\curl)$ Riesz maps~\cite{arnold1997, arnold2000}.

However, the proof that this element pair satisfies the kernel decomposition property
depends on the pressure space being piecewise constant, and does not generalise to either Taylor--Hood or
conforming divergence-free finite elements.  In fact, the same choice of
space decomposition applied to a barycentrically refined mesh does not result
in a $\gamma$-robust smoother for the Scott--Vogelius element considered here.

To find a space decomposition that decomposes the kernel of the divergence
operator for the Scott--Vogelius element, we consider the following Hilbert
complexes in two
\begin{equation}
  \label{eq:stokes-complex-2d_2}
  \mathbb{R} \xrightarrow{\operatorname{id}} H^2
  \xrightarrow{\curl} H^1 \xrightarrow{\div} L^2 \xrightarrow{\operatorname{null}} 0
\end{equation}
and three dimensions
\begin{equation}
  \label{eq:stokes-complex-3d_3}
  \mathbb{R} \xrightarrow{\operatorname{id}} H^2
  \xrightarrow{\grad} H^1(\curl) \xrightarrow{\curl} H^1 \xrightarrow{\div} L^2 \xrightarrow{\operatorname{null}} 0,
\end{equation}
where $H^1(\curl) = \{u \in (H^1)^3, \curl u \in (H^1)^3 \}$.
These sequences are exact on simply connected domains, which
implies that every field $u$ in the kernel of the divergence can be represented as the $\curl$ of a potential $\Phi$.
If we are given a discrete subcomplex of the form
\begin{equation}
 \ldots \rightarrow \Sigma_h \xrightarrow{\curl} V_h \xrightarrow{\div} Q_h \xrightarrow{\operatorname{null}} 0,
\end{equation}
then for a divergence-free discrete vector-field $u_h\in V_h$, we can write it as the $\curl$ of a potential $\Phi_h\in\Sigma_h$.
Writing $\Phi_h=\sum_j \Phi_j $ in terms of basis functions of $\Sigma_h$, we
can then define a divergence-free decomposition of $u_h$ as $u_h = \sum_j \nabla\times \Phi_j$.
Hence, a space decomposition such that $\nabla \times \Phi_j$ is contained in
some $V_i$ for any basis function $\Phi_j$ decomposes the kernel and satisfies
\eqref{eqn:kernel-decomposition}.

\begin{figure}[htbp]
  \centering
  \begin{tikzpicture}[scale=4,
    normaldofs/.style={postaction=decorate, decoration={markings, mark=at position 0.5 with {
            \draw (0,-5pt) -- ++(0, 10pt);
    }}},
    pointdofs/.style={postaction=decorate, decoration={markings, mark=at position 0.5 with {
                \draw [fill] circle (0.12);
    }}}
    ]
    \coordinate (0) at (0, 0) {};
    \coordinate (1) at (1, 0) {};
    \coordinate (2) at ($(0) + (60:1)$) {};
    \coordinate (bc) at (barycentric cs:0=1,1=1,2=1) {};

    \draw[normaldofs] (0) to (1);
    \draw[normaldofs] (1) to (2);
    \draw[normaldofs] (2) to (0);

    \foreach \i in {0, 1, 2} {
        \draw (\i) circle (0.06);
        \draw[fill] (\i) circle (0.03);
    }

    \foreach \i in {0, 1, 2} {
        \draw[dashed] (bc) -- (\i);
    }
    \begin{scope}[shift={(1.5, 0)}]

        \coordinate (0) at (0, 0) {};
        \coordinate (1) at (1, 0) {};
        \coordinate (2) at ($(0) + (60:1)$) {};
        \coordinate (bc) at (barycentric cs:0=1,1=1,2=1) {};

        \draw[pointdofs] (0) to (1);
        \draw[pointdofs] (1) to (2);
        \draw[pointdofs] (2) to (0);
        \foreach \i in {0, 1, 2} {
            \draw[pointdofs] (bc) -- (\i);
        }

        \foreach \i in {0, 1, 2} {
            \draw[fill] (\i) circle (0.03);
        }

        \draw[fill] (bc) circle (0.03);
    \end{scope}
    \begin{scope}[shift={(3, 0)}]
        \coordinate (0) at (0, 0) {};
        \coordinate (1) at (1, 0) {};
        \coordinate (2) at ($(0) + (60:1)$) {};
        \coordinate (bc) at (barycentric cs:0=1,1=1,2=1) {};
        \draw[line join=miter] (0) -- (1) -- (2) -- cycle;
        \draw[line join=miter] (0) -- (bc) -- (1);
        \draw[line join=miter] (2) -- (bc);

        \coordinate (bc0) at (barycentric cs:0=1,1=1,bc=1) {};
        \coordinate (bc1) at (barycentric cs:1=1,2=1,bc=1) {};
        \coordinate (bc2) at (barycentric cs:0=1,2=1,bc=1) {};

        \foreach \i in {0, 1, 2} {
            \draw[fill] ($(bc\i)-(\i*120 + 0:0.07)$) circle (0.03);
            \draw[fill] ($(bc\i)+(\i*120 + 0:0.07)$) circle (0.03);
            \draw[fill] ($(bc\i)+(\i*120 + 90:0.07)$) circle (0.03);
        }
    \end{scope}
\end{tikzpicture}
  \caption{\label{fig:2d-complex}2D exact Stokes complex.}
\end{figure}
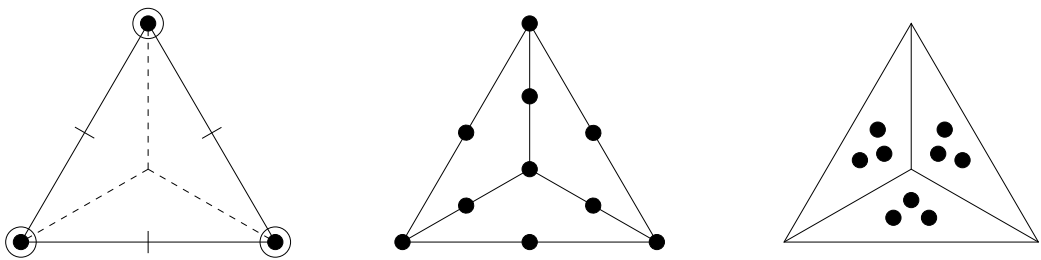

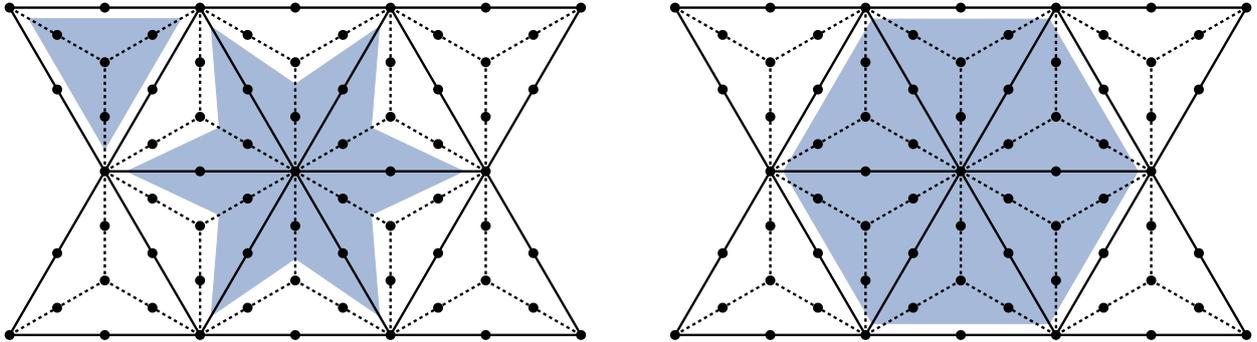
\begin{figure}[!htbp]
    \centering
    \begin{tikzpicture}[scale=4,
  bcdof/.style={postaction={decorate,decoration={markings,
        mark=at position 0.5 with {\draw[solid, fill, black] circle
          (0.12);}}}},
  dof/.style={postaction={decorate,decoration={markings,
        mark=at position 0.5 with {\draw[solid, fill, black] circle (0.12);}}}
  }]

  \foreach \i in {0, 1, 2, 3} {
      \coordinate(0-\i) at ($(\i*1.5, 0)$);
      \coordinate(1-\i) at ($(0-\i) + (60:1.5)$);
      \coordinate(2-\i) at ($(1-\i) + (120:1.5)$);

    \foreach \j in {0, 2} {
      \draw[line width=0.75mm, solid, black, fill] (\j-\i) circle (0.03);
    }
  }
  \foreach \i in {0, 1, 2} {
    \draw[line width=0.75mm, solid, black, fill] (1-\i) circle (0.03);
  }

  \foreach \i/\j in {0/1, 1/2, 2/3} {
    \foreach \k in {0, 2} {
      \draw[dof, line join=miter, line width=0.75mm] (\k-\i) -- (\k-\j);
    }
    \draw[dof, line join=miter, line width=0.75mm] (0-\i) -- (1-\i);
    \draw[dof, line join=miter, line width=0.75mm] (1-\i) -- (2-\i);
    \draw[dof, line join=miter, line width=0.75mm] (0-\j) -- (1-\i);
    \draw[dof, line join=miter, line width=0.75mm] (1-\i) -- (2-\j);
  }

  \foreach \i/\j in {0/1, 1/2} {
    \draw[dof, line join=miter, line width=0.75mm] (1-\i) -- (1-\j);
  }
  
  \foreach \i/\j in {0/1, 1/2, 2/3} {
      \coordinate (bc-0-\i-\j) at (barycentric cs:0-\i=1,0-\j=1,1-\i=1);
    \draw[bcdof, line join=miter, dashed, line width=0.75mm] (0-\i) -- (bc-0-\i-\j);
    \draw[bcdof, line join=miter, dashed, line width=0.75mm] (1-\i) -- (bc-0-\i-\j);
    \draw[bcdof, line join=miter, dashed, line width=0.75mm] (0-\j) -- (bc-0-\i-\j);

    \coordinate (bc-2-\i-\j) at (barycentric cs:1-\i=1,2-\i=1,2-\j=1);
    \draw[bcdof, line join=miter, dashed, line width=0.75mm] (1-\i) -- (bc-2-\i-\j);
    \draw[bcdof, line join=miter, dashed, line width=0.75mm] (2-\i) -- (bc-2-\i-\j);
    \draw[bcdof, line join=miter, dashed, line width=0.75mm] (2-\j) -- (bc-2-\i-\j);

    \foreach \k in {0, 2} {
      \draw[line width=0.75mm, solid, black, fill] (bc-\k-\i-\j) circle (0.03);
    }
  }
  \foreach \i/\j in {0/1, 1/2} { 
    \coordinate (bc-1-\i-\j)  at (barycentric cs:1-\i=1,0-\j=1,1-\j=1);
    \draw[bcdof, line join=miter, dashed, line width=0.75mm] (1-\i) -- (bc-1-\i-\j);
    \draw[bcdof, line join=miter, dashed, line width=0.75mm] (0-\j) -- (bc-1-\i-\j);
    \draw[bcdof, line join=miter, dashed, line width=0.75mm] (1-\j) -- (bc-1-\i-\j);

    \coordinate (bc-3-\i-\j)  at (barycentric cs:1-\i=1,1-\j=1,2-\j=1);
    \draw[bcdof, line join=miter, dashed, line width=0.75mm] (1-\i) -- (bc-3-\i-\j);
    \draw[bcdof, line join=miter, dashed, line width=0.75mm] (1-\j) -- (bc-3-\i-\j);
    \draw[bcdof, line join=miter, dashed, line width=0.75mm] (2-\j) -- (bc-3-\i-\j);
    \foreach \k in {1, 3} {
      \draw[line width=0.75mm, solid, black, fill] (bc-\k-\i-\j) circle (0.03);
    }
  }

  \begin{scope}[on background layer]
      \draw[blue!50, fill=blue!50] ([shift={(0:5pt)}]1-0) --
      ([shift={(30:5pt)}] bc-1-0-1) --
      ([shift={(60:5pt)}]0-1) --
      ([shift={(90:5pt)}] bc-0-1-2) --
      ([shift={(120:5pt)}]0-2) --
      ([shift={(150:5pt)}] bc-1-1-2) --
      ([shift={(180:5pt)}]1-2) --
      ([shift={(210:5pt)}] bc-3-1-2) --
      ([shift={(240:5pt)}]2-2) --
      ([shift={(270:5pt)}] bc-2-1-2) --
      ([shift={(300:5pt)}]2-1) --
      ([shift={(330:5pt)}] bc-3-0-1) --
      cycle;

      \draw[blue!50, fill=blue!50] ([shift={(330:5pt)}]2-0) --
      ([shift={(210:5pt)}]2-1) --
      ([shift={(90:5pt)}]1-0) --
      cycle;
  \end{scope}
\end{tikzpicture}
    \hspace{0.05\textwidth}
    \begin{tikzpicture}[scale=4,
  bcdof/.style={postaction={decorate,decoration={markings,
        mark=at position 0.5 with {\draw[solid, fill, black] circle
          (0.12);}}}},
  dof/.style={postaction={decorate,decoration={markings,
        mark=at position 0.5 with {\draw[solid, fill, black] circle (0.12);}}}
  }]

  \foreach \i in {0, 1, 2, 3} {
    \coordinate (0-\i) at ($(\i*1.5, 0)$);
    \coordinate (1-\i) at ($(0-\i) + (60:1.5)$);
    \coordinate (2-\i) at ($(1-\i) + (120:1.5)$);

    \foreach \j in {0, 2} {
      \draw[line width=0.75mm, solid, fill, black] (\j-\i) circle (0.03);
    }
  }
  \foreach \i in {0, 1, 2} {
    \draw[line width=0.75mm, solid, fill, black] (1-\i) circle (0.03);
  }

  \foreach \i/\j in {0/1, 1/2, 2/3} {
    \foreach \k in {0, 2} {
      \draw[dof, line join=miter, line width=0.75mm] (\k-\i) -- (\k-\j);
    }
    \draw[dof, line join=miter, line width=0.75mm] (0-\i) -- (1-\i);
    \draw[dof, line join=miter, line width=0.75mm] (1-\i) -- (2-\i);
    \draw[dof, line join=miter, line width=0.75mm] (0-\j) -- (1-\i);
    \draw[dof, line join=miter, line width=0.75mm] (1-\i) -- (2-\j);
  }

  \foreach \i/\j in {0/1, 1/2} {
    \draw[dof, line join=miter, line width=0.75mm] (1-\i) -- (1-\j);
  }
  
  \foreach \i/\j in {0/1, 1/2, 2/3} {
    \coordinate (bc-0-\i-\j) at (barycentric cs:0-\i=1,0-\j=1,1-\i=1);
    \draw[bcdof, line join=miter, dashed, line width=0.75mm] (0-\i) -- (bc-0-\i-\j);
    \draw[bcdof, line join=miter, dashed, line width=0.75mm] (1-\i) -- (bc-0-\i-\j);
    \draw[bcdof, line join=miter, dashed, line width=0.75mm] (0-\j) -- (bc-0-\i-\j);

    \coordinate (bc-2-\i-\j) at (barycentric cs:1-\i=1,2-\i=1,2-\j=1);
    \draw[bcdof, line join=miter, dashed, line width=0.75mm] (1-\i) -- (bc-2-\i-\j);
    \draw[bcdof, line join=miter, dashed, line width=0.75mm] (2-\i) -- (bc-2-\i-\j);
    \draw[bcdof, line join=miter, dashed, line width=0.75mm] (2-\j) -- (bc-2-\i-\j);

    \foreach \k in {0, 2} {
      \draw[line width=0.75mm, solid, fill, black] (bc-\k-\i-\j) circle (0.03);
    }
  }
  \foreach \i/\j in {0/1, 1/2} { 
    \coordinate (bc-1-\i-\j)  at (barycentric cs:1-\i=1,0-\j=1,1-\j=1);
    \draw[bcdof, line join=miter, dashed, line width=0.75mm] (1-\i) -- (bc-1-\i-\j);
    \draw[bcdof, line join=miter, dashed, line width=0.75mm] (0-\j) -- (bc-1-\i-\j);
    \draw[bcdof, line join=miter, dashed, line width=0.75mm] (1-\j) -- (bc-1-\i-\j);

    \coordinate (bc-3-\i-\j)  at (barycentric cs:1-\i=1,1-\j=1,2-\j=1);
    \draw[bcdof, line join=miter, dashed, line width=0.75mm] (1-\i) -- (bc-3-\i-\j);
    \draw[bcdof, line join=miter, dashed, line width=0.75mm] (1-\j) -- (bc-3-\i-\j);
    \draw[bcdof, line join=miter, dashed, line width=0.75mm] (2-\j) -- (bc-3-\i-\j);
    \foreach \k in {1, 3} {
      \draw[line width=0.75mm, solid, fill, black] (bc-\k-\i-\j) circle (0.03);
    }
  }




  \begin{scope}[on background layer]
      \draw[blue!50, fill=blue!50] ([shift={(0:3pt)}]1-0) --
  ([shift={(60:3pt)}]0-1) --
  ([shift={(120:3pt)}]0-2) --
  ([shift={(180:3pt)}]1-2) --
  ([shift={(240:3pt)}]2-2) --
  ([shift={(300:3pt)}]2-1) --
  cycle;
  \end{scope}
\end{tikzpicture}
    \caption{Left: The star patch applied to the barycentrically refined mesh does not yield a robust relaxation method for the Scott--Vogelius element. Right: The $\MacroStar$ patch satisfies the kernel decomposition property \eqref{eqn:kernel-decomposition} for the $[\Pk]^d\mathrm{-}\Pkminusdisc$ element for $k\ge d$ (here shown for $k=d=2$).\label{fig:macrostar}}
\end{figure}

In two dimensions, the Scott--Vogelius velocity and pressure spaces $V_h$ and
$Q_h$ form an exact sequence with $\Sigma_h$ chosen as the HCT finite element space \cite[\S 6.1]{ciarlet1978}.
The three elements are displayed in Figure~\ref{fig:2d-complex}.
For a given vertex $v_i$ in the macro mesh $\mesh_h$, we define the $\MacroStar(v_i)$ of the vertex as the union of all macro elements touching the vertex.
We then see that for every basis function $\Phi_j$ there exists a vertex $v_i$ such that $\mathrm{supp}(\Phi_j)\subset \MacroStar(v_i)$.
Hence also $\mathrm{supp}(\nabla\times \Phi_j)\subset\MacroStar(v_i)$ and if we define
\begin{equation}
    V_i = \{v\in V_h : \mathrm{supp}(v)\subset\MacroStar(v_i)\}
\end{equation}
then these subspaces decompose the kernel.
The $\MacroStar$ is shown in Figure~\ref{fig:macrostar}.

A similar argument can be made in three dimensions. Here the existence of an
exact sequence is given by the recent work of Fu, Guzm\'an, and Neilan~\cite{fu_exact_2020}.

We note that the proof for robustness of the multigrid scheme in
\cite{schoberl1999b} actually has a stricter requirement than simply the kernel
decomposition property \eqref{eqn:kernel-decomposition}. The kernel
decomposition must be \emph{stable}: given a function $u\in \mathcal{N}_h$ with
$u=\sum u_i$ and $u_i \in V_i \cap \mathcal{N}_h$, one needs to be able to
estimate the norm of the $u_i$ in terms of $u$. In general this estimate does
not follow purely from the exactness of the discrete sequence, but it was shown
in \cite{fmsw_robust_2020} that one can use the existence of a particular Fortin
operator to obtain the required bounds.

\subsection{Relaxation in the presence of stabilization}
\label{sec:stabilization}
It is well known that straightforward Galerkin discretizations of advection-dominated problems are
oscillatory~\cite{brooks1982,turek1999,quarteroni2008,elman2014}.
Several approaches have been developed to address these issues, for example by adding a small amount of artificial viscosity as in the case of Streamline Upwind/Petrov Galerkin (SUPG) or Galerkin Least Squares (GLS) or by enriching the space with bubble functions.
SUPG stabilization was used in previous works
on Reynolds-robust preconditioners~\cite{benzi2006,farrell2018b}.
As both of these works consider piecewise constant pressures, the pressure
gradient on each cell vanishes and hence only the top-left block is modified when adding SUPG stabilization.
For the Scott--Vogelius element this is not the case, and SUPG stabilization results in a modification of the top-right block, introducing additional coupling between the velocity and the pressure~\cite[p.~1249]{burman_continuous_2006},~\cite[p.~250]{gelhard_stabilized_2005}.
This makes the nonlinear problem much more difficult to solve. In fact, in
numerical experiments we observe that the outer Newton solver already fails to converge
at $\Re = 50$ for a two dimensional lid-driven cavity, even with the
use of a direct solver.
Furthermore, the modified top-right block needs to be taken into account when
adding the augmented Lagrangian term and in the multigrid scheme for the
top-left block.

In 1976 Douglas \& Dupont~\cite{douglas1976} suggested an interior penalty scheme that penalises a jump of the derivative across facets:
\begin{equation}\label{eqn:burman-stab}
    S(u, v) = \sum_{K\in\hat\mesh_h} \frac{1}{2} \int_{\partial K}  \delta
    h_{\partial K}^2 \dgjmpls\nabla u\dgjmprs : \dgjmpls\nabla v\dgjmprs \dr s,
\end{equation}
where $\dgjmpls\nabla u\dgjmprs$ denotes the jump of the gradient, $h_{\partial K}$ is a function giving the facet size, and $\delta$ is a free parameter.
The term vanishes when the velocity is $C^1$-continuous.
This scheme has received renewed attention and it was shown in~\cite{burman_edge_2004} that it successfully stabilizes advection-dominated problems
and has subsequently been used to stabilize the Stokes~\cite{burman_edge_2006} and the Oseen equations~\cite{burman_continuous_2006, burman_stabilized_2008}.
\revnew{In addition to not introducing any additional coupling of the pressure and the velocity, an advantage of this scheme is that it is adjoint consistent.
This means that for low Reynolds number one can prove that the order of convergence in the $L^2$ norm is preserved~\cite[Remark 12]{burman_stabilized_2008}.
In fact, for the problem considered in Section~\ref{sec:svmms} we observe optimal convergence even at high Reynolds number, although this is not in general guaranteed.}

We now consider the effect of adding~\eqref{eqn:burman-stab} to the top-left
block on the multigrid scheme.
Since $S$ vanishes for functions that have continuous gradients, we have added another bilinear form to our system that has a nontrivial kernel consisting of $C^1$ vector fields.
As the weight $\delta h_{\partial K}^2$ is small, the impact is not as significant as that of the grad-div term, but for very high Reynolds number or coarse meshes, we still observe reduced performance of the multigrid scheme.
As discussed in the previous section, we know that for the smoother to be robust the space decomposition must provide a decomposition of the kernel.
In two dimensions, this is satisfied if $k\ge 3$, as the $\MacroStar$ around vertices then captures the support of the HCT element.
In three dimensions the lowest degree (that the authors are aware of) for a local basis for $C^1$ vector fields on barycentrically refined meshes is $k=5$, see~\cite{alfeld_1984,fasshauer_multivariate_2014}.

This argument is heuristic, and a full analysis of problems with two different
singular terms is out of the scope of this work. Nevertheless, in the numerical
experiments we will see that the scheme is noticeably more robust for $k=3$ in
two dimensions and $k=5$ in three dimensions in the presence of the
stabilization \eqref{eqn:burman-stab}.

\begin{remark}
  We conjecture that similar robustness properties will carry over to the case
  of Powell--Sabin splits. In two dimensions, there is a local
  $C^1$-conforming quadratic basis on Powell--Sabin splits
  \cite{powell_sabin_1977}, and in three dimensions there is a local
  $C^1$-conforming cubic basis \cite{worsey_farin_1987}. We
  therefore expect that the choice $k \ge d$ will provide more robust iteration counts in the presence of
  stabilization than the minimal $k = d - 1$ required for inf-sup stability on meshes with this
  macro structure.
\end{remark}

\subsection{Prolongation}
The second key ingredient of the multigrid scheme is a robust prolongation operator.
To keep notation simple we consider the case of a two-level scheme and denote the coarse-grid function space by $V_H$ and the fine-grid space by $V_h$.
We denote the standard prolongation operator induced by the finite element
interpolation operator in $V_h$ by $P_H:V_H\to V_h$.
Let $u_H\in V_H$ be a coarse-grid function.
For the multigrid scheme to be $\gamma$-robust, it was shown by Sch\"oberl \cite{schoberl1999b} that the prolongation operator
must satisfy
\begin{equation}
    \| P_H u_H\|_{A_{h,\gamma}} \le C \|u_H\|_{A_{H,\gamma}}
\end{equation}
with a constant $C$ independent of $\gamma$.
Calculating these norms for a divergence-free function $u_H\in V_H$, we observe
that
\begin{equation}
    \begin{aligned}
        \|u_H\|_{A_{H,\gamma}}^2    &= \|\nabla u_H\|_{L^2}^2 + \gamma \|\xunderbrace{\nabla\cdot u_H}_{=0}\|_{L^2}^2,\\
        \|P_Hu_H\|_{A_{h,\gamma}}^2 &= \|\nabla (P_Hu_H)\|_{L^2}^2 + \gamma \| \nabla\cdot (P_Hu_H)\|_{L^2}^2.
    \end{aligned}
\end{equation}
However, since the multigrid hierarchy considered here is non-nested, the interpolation is not exact and $u_H$ being divergence-free does not necessarily imply the same property for $P_Hu_H$.
\begin{remark}
   We note that this scenario differs from the situation in \cite{schoberl1999b, benzi2006}.
   There the function spaces are nested, but a discretely divergence-free
   function on the coarse-grid may not be discretely divergence-free on the fine
   grid, as $Q_h$ is larger than $Q_H$.
\end{remark}

\begin{figure}[htbp]
  \centering
  \ExplSyntaxOn
\newcommand{\fpeval}[1]{\fp_eval:n{#1}}
\ExplSyntaxOff

\newcommand{\convexpath}[2]{
  [
  create hullcoords/.code={
    \global\edef\namelist{#1}
    \foreach [count=\counter] \nodename in \namelist {
      \global\edef\numberofnodes{\counter}
      \coordinate (hullcoord\counter) at (\nodename);
    }
    \coordinate (hullcoord0) at (hullcoord\numberofnodes);
    \pgfmathtruncatemacro\lastnumber{\numberofnodes+1}
    \coordinate (hullcoord\lastnumber) at (hullcoord1);
  },
  create hullcoords
  ]
  ($(hullcoord1)!#2!-90:(hullcoord0)$)
  \foreach [
  evaluate=\currentnode as \previousnode using \currentnode-1,
  evaluate=\currentnode as \nextnode using \currentnode+1
  ] \currentnode in {1,...,\numberofnodes} {
    let \p1 = ($(hullcoord\currentnode) - (hullcoord\previousnode)$),
    \n1 = {atan2(\y1,\x1) + 90},
    \p2 = ($(hullcoord\nextnode) - (hullcoord\currentnode)$),
    \n2 = {atan2(\y2,\x2) + 90},
    \n{delta} = {Mod(\n2-\n1,360) - 360}
    in
    {arc [start angle=\n1, delta angle=\n{delta}, radius=#2]}
    -- ($(hullcoord\nextnode)!#2!-90:(hullcoord\currentnode)$)
  }
}

\newcommand{\DrawIterated}[8]
{ \ifnum#1=1%
  \draw[line join=round, line width=#8, dof] (#2,#3) -- (#4,#5) -- (#6,#7) -- cycle;
  \else
  \draw[line join=round, line width=#8, dof] (#2,#3) -- (#4,#5) -- (#6,#7) -- cycle;
  \DrawIterated{\numexpr#1-1\relax}%
  {\fpeval{(#2 + #4)/2}}%
  {\fpeval{(#3 + #5)/2}}%
  {\fpeval{(#4 + #6)/2}}%
  {\fpeval{(#5 + #7)/2}}%
  {\fpeval{(#2 + #6)/2}}%
  {\fpeval{(#3 + #7)/2}}

  \DrawIterated{\numexpr#1-1\relax}%
  {#2}%
  {#3}%
  {\fpeval{(#2 + #4)/2}}%
  {\fpeval{(#3 + #5)/2}}%
  {\fpeval{(#2 + #6)/2}}%
  {\fpeval{(#3 + #7)/2}}%
  {#8}

  \DrawIterated{\numexpr#1-1\relax}%
  {\fpeval{(#2 + #4)/2}}%
  {\fpeval{(#3 + #5)/2}}%
  {#4}%
  {#5}%
  {\fpeval{(#4 + #6)/2}}%
  {\fpeval{(#5 + #7)/2}}%
  {#8}

  \DrawIterated{\numexpr#1-1\relax}%
  {\fpeval{(#2 + #6)/2}}%
  {\fpeval{(#3 + #7)/2}}
  {\fpeval{(#4 + #6)/2}}%
  {\fpeval{(#5 + #7)/2}}%
  {#6}%
  {#7}%
  {#8}
  
  \fi
}

\newcommand{\DrawIteratedBary}[7]
{ \ifnum#1=1%
  \draw[line join=round, line width=0.2mm, black, dof] (#2,#3) -- (#4,#5);
  \draw[line join=round, line width=0.2mm, black, dof] (#4,#5) -- (#6,#7);
  \draw[line join=round, line width=0.2mm, black, dof] (#6,#7) -- (#2,#3);
  \coordinate (0) at (#2,#3);
  \coordinate (1) at (#4,#5);
  \coordinate (2) at (#6,#7);
  \coordinate (bc) at (barycentric cs:0=1,1=1,2=1);
  \draw[line join=round, black, line width=0.2mm, dashed, dof] (0) -- (bc);
  \draw[line join=round, black, line width=0.2mm, dashed, dof] (1) -- (bc);
  \draw[line join=round, black, line width=0.2mm, dashed, dof] (2) -- (bc);

  \else
  \draw[line join=round, black] (#2,#3) -- (#4,#5) -- (#6,#7) -- (#2, #3);
  \DrawIteratedBary{\numexpr#1-1\relax}%
  {\fpeval{(#2 + #4)/2}}%
  {\fpeval{(#3 + #5)/2}}%
  {\fpeval{(#4 + #6)/2}}%
  {\fpeval{(#5 + #7)/2}}%
  {\fpeval{(#2 + #6)/2}}%
  {\fpeval{(#3 + #7)/2}}

  \DrawIteratedBary{\numexpr#1-1\relax}%
  {#2}%
  {#3}%
  {\fpeval{(#2 + #4)/2}}%
  {\fpeval{(#3 + #5)/2}}%
  {\fpeval{(#2 + #6)/2}}%
  {\fpeval{(#3 + #7)/2}}

  \DrawIteratedBary{\numexpr#1-1\relax}%
  {\fpeval{(#2 + #4)/2}}%
  {\fpeval{(#3 + #5)/2}}%
  {#4}%
  {#5}%
  {\fpeval{(#4 + #6)/2}}%
  {\fpeval{(#5 + #7)/2}}%

  \DrawIteratedBary{\numexpr#1-1\relax}%
  {\fpeval{(#2 + #6)/2}}%
  {\fpeval{(#3 + #7)/2}}
  {\fpeval{(#4 + #6)/2}}%
  {\fpeval{(#5 + #7)/2}}%
  {#6}%
  {#7}
  
  \fi
}

\begin{tikzpicture}[scale=4,
  dof/.style={
  postaction={decorate,decoration={markings,
      mark=at position 0.0 with {\draw[solid, fill, black] circle (0.075cm);},
      mark=at position 0.5 with {\draw[solid, fill, black] circle (0.075cm);},
      mark=at position 1.0 with {\draw[solid, fill, black] circle (0.075cm);}}}}]
  \DrawIteratedBary{2}{1}{0}{0.5}{\fpeval{sqrt(3)/2}}{1.5}{\fpeval{sqrt(3)/2}};
  \DrawIteratedBary{2}{1}{0}{2}{0}{1.5}{\fpeval{sqrt(3)/2}};
  \DrawIteratedBary{2}{2}{0}{1.5}{\fpeval{sqrt(3)/2}}{2.5}{\fpeval{sqrt(3)/2}};

  \DrawIteratedBary{2}{0.5}{\fpeval{sqrt(3)/2}}{1.5}{\fpeval{sqrt(3)/2}}{1}{\fpeval{sqrt(3)}};
  \DrawIteratedBary{2}{1.5}{\fpeval{sqrt(3)/2}}{2}{\fpeval{sqrt(3)}}{1}{\fpeval{sqrt(3)}};
  \DrawIteratedBary{2}{1.5}{\fpeval{sqrt(3)/2}}{2.5}{\fpeval{sqrt(3)/2}}{2}{\fpeval{sqrt(3)}};

  \DrawIterated{1}{1}{0}{0.5}{\fpeval{sqrt(3)/2}}{1.5}{\fpeval{sqrt(3)/2}}{0.6mm};
  \DrawIterated{1}{1}{0}{2}{0}{1.5}{\fpeval{sqrt(3)/2}}{0.6mm};
  \DrawIterated{1}{2}{0}{1.5}{\fpeval{sqrt(3)/2}}{2.5}{\fpeval{sqrt(3)/2}}{0.6mm};

  \DrawIterated{1}{0.5}{\fpeval{sqrt(3)/2}}{1.5}{\fpeval{sqrt(3)/2}}{1}{\fpeval{sqrt(3)}}{0.6mm};
  \DrawIterated{1}{1.5}{\fpeval{sqrt(3)/2}}{2}{\fpeval{sqrt(3)}}{1}{\fpeval{sqrt(3)}}{0.6mm};
  \DrawIterated{1}{1.5}{\fpeval{sqrt(3)/2}}{2.5}{\fpeval{sqrt(3)/2}}{2}{\fpeval{sqrt(3)}}{0.6mm};

  \begin{scope}[on background layer]


    \coordinate (d) at ($(1, 0) + (90:0.075)$);
    \coordinate (e) at ($(0.5, {sqrt(3)/2}) + (-30:0.075)$);
    \coordinate (f) at ($(1.5, {sqrt(3)/2}) + (210:0.075)$);

    \draw[blue!50, fill=blue!50] \convexpath{d,e,f}{0pt};

    \begin{scope}[shift={(1, 0)}]


      \coordinate (d) at ($(1, 0) + (90:0.075)$); \coordinate (e) at
      ($(0.5, {sqrt(3)/2}) + (-30:0.075)$); \coordinate (f) at
      ($(1.5, {sqrt(3)/2}) + (210:0.075)$);

      \draw[blue!50, fill=blue!50] \convexpath{d,e,f}{0pt};
    \end{scope}


    \begin{scope}[yscale=-1, shift={(0, -{sqrt(3)})}]




      \begin{scope}[shift={(1, 0)}]
        \coordinate (a) at ($(0, 0) + (30:0.075)$); \coordinate (b) at
        ($(0.5, {sqrt(3)/2}) + (270:0.075)$); \coordinate (c) at
        ($(1, 0) + (150:0.075)$);

        \draw[blue!50, fill=blue!50] \convexpath{a,b,c}{0pt};


      \end{scope}
      \begin{scope}[shift={(2, 0)}]


      \end{scope}

    \end{scope}
  \end{scope}

\end{tikzpicture}
  \caption{The robust prolongation operator solves local Stokes problems on each coarse macro cell.}
  \label{fig:prolongation-patch-2d}
\end{figure}
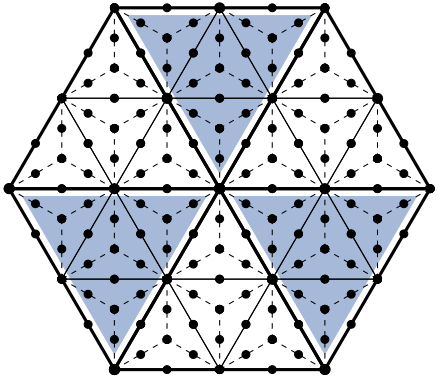

Inspecting the mesh hierarchy in Figure~\ref{fig:bary-hierarchy}, we notice that
the interpolation is exact along the edges of the coarse-grid macro mesh.  In
turn, this means that the flux across these edges is preserved exactly and hence
the interpolated vector-field is divergence-free with respect to pressure
functions that are piecewise constant on the macro mesh (i.e.~before barycentric
refinement):
\begin{equation}
    \tilde Q_H = \{ q \in L^2(\Omega): q\vert_{K} \equiv \text{const} \ \forall K \in \mesh_h\}.
\end{equation}
For a robust prolongation, we therefore only need to modify the degrees of freedom within a coarse-grid macro cell to remove the divergence within the cell that was created by the interpolation.
To this end we define the subspace $\tilde V_h\subset V_h$ of functions that
vanish on the boundaries of macro cells
\begin{equation}
    \tilde V_h = \{v \in V_h : v =0 \text{ on } \partial K \text{ for all } K \in \mesh_h\}.
\end{equation}
We then solve for $\tilde u_h\in \tilde V_h$ such that
\begin{equation}\label{eqn:prolongation-local-problem}
    \nu (\eps{\tilde u_h}, \eps{\tilde v_h}) + \gamma (\nabla\cdot \tilde u_h, \nabla\cdot\tilde v_h)= \gamma (\Pi_{\tilde Q_h} (\nabla\cdot (P_H u_H)), \Pi_{\tilde Q_h} (\nabla\cdot \tilde v_h)) \quad \text{for all } \tilde v_h \in \tilde V_h.
\end{equation}
It was shown in \cite{fmsw_robust_2020} that then the modified prolongation given by
\begin{equation}
    \tilde P_H u_H = P_Hu_H - \tilde u_h
\end{equation}
is continuous in the energy norm with a continuity constant independent of
$\gamma$.
We emphasize that due to the nature of the space $\tilde V_h$,
\eqref{eqn:prolongation-local-problem} decouples into many small, independent
solves and hence can be solved efficiently, see Figure~\ref{fig:prolongation-patch-2d}.
In addition, since the problem solved on each coarse macro cell does not vary
through the nonlinear iteration, the small matrices can be assembled and
factorized once in an initialization step.

\section{Numerical examples} \label{sec:examples}
\subsection{Details of the algorithm}

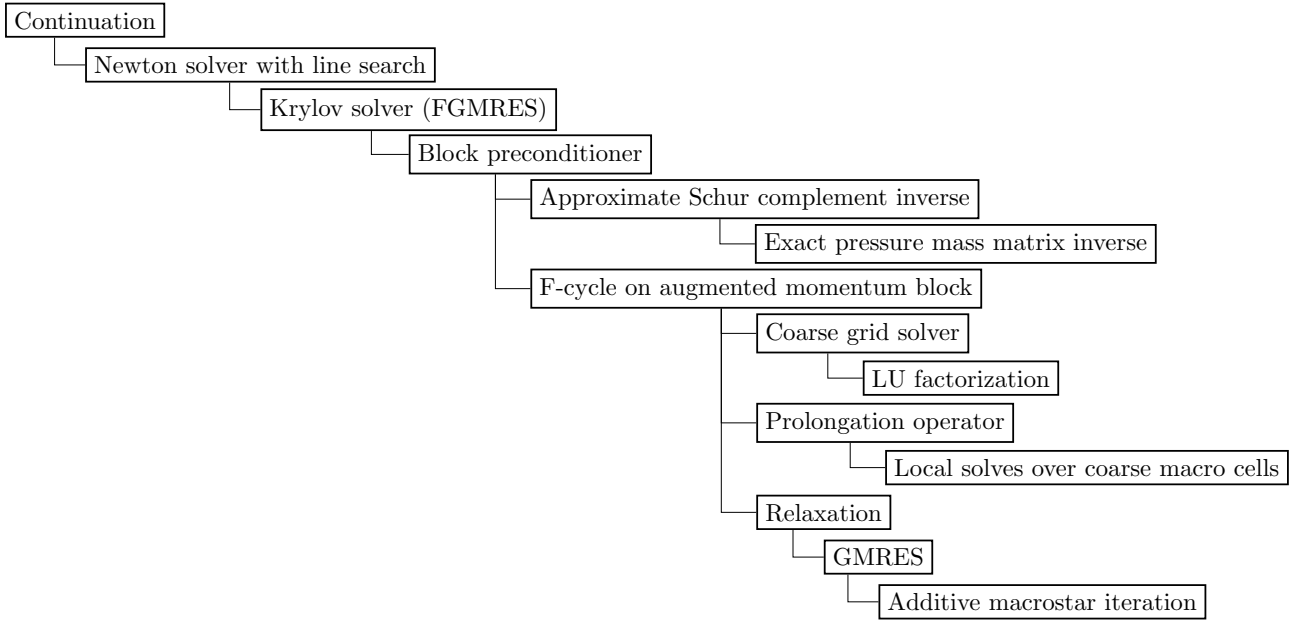
\begin{figure}[htbp]
  \centering
  \resizebox{\textwidth}{!}{\begin{tikzpicture}[%
  every node/.style={draw=black, thick, anchor=west},
  grow via three points={one child at (0.0,-0.7) and
  two children at (0.0,-0.7) and (0.0,-1.4)},
  edge from parent path={(\tikzparentnode.210) |- (\tikzchildnode.west)}]
  \node {Continuation}
    child { node {Newton solver with line search}
      child { node {Krylov solver (FGMRES)}
        child { node {Block preconditioner}
          child { node {Approximate Schur complement inverse}
              child { node {Exact pressure mass matrix inverse}}
          }
          child [missing] {}
          child { node {F-cycle on augmented momentum block}
              child { node {Coarse grid solver}
                child { node {LU factorization}}
              }
              child [missing] {}
              child { node {Prolongation operator}
                child { node {Local solves over coarse macro cells}}
              }
              child [missing] {}
              child { node {Relaxation}
                child { node {GMRES}
                  child { node {Additive macrostar iteration}}
                }
              }
          }
        }
      }
    };
\end{tikzpicture}}
  \caption{An outline of the algorithm for solving \eqref{eq:ns}.}
  \label{fig:solver}
\end{figure}

A graphical depiction of the solver algorithm is given in Figure \ref{fig:solver}.
The essential structure is the same as in \cite{farrell2018b}, with different
multigrid components for $A_\gamma$.
The code is implemented in Firedrake
\cite{rathgeber2016} using PETSc \cite{balay2019} and PCPATCH \cite{farrell2019c}.
Since the convergence behavior of Newton's method is not Reynolds-robust, we employ continuation
in the Reynolds number to ensure its convergence.
A flexible Krylov variant is required as we apply GMRES inside the multigrid relaxation,
and hence we use flexible GMRES \cite{saad1993} as the outermost Krylov solver.
We use the full block factorization
preconditioner
\begin{equation}
P^{-1} =
\left( \begin{array}{cc}
I   & -\tilde{A}_\gamma^{-1} B^\top \\
0 & I \\
\end{array} \right)
\left( \begin{array}{cc}
\tilde{A}_\gamma^{-1}  & 0 \\
0 & \tilde{S}^{-1} \\
\end{array} \right)
\left( \begin{array}{cc}
I   & 0 \\
-B\tilde{A}_\gamma^{-1} & I \\
\end{array} \right)
\end{equation}
with the scaled inverse of the (block diagonal) pressure mass matrix as $\tilde{S}^{-1}$,
and one full multigrid cycle of the algorithm described in Section \ref{sec:multigrid}
as $\tilde{A}_{\gamma}^{-1}$. Each relaxation sweep conducts 6 (in 2D) or 10 (in 3D) GMRES iterations
preconditioned by the additive macrostar iteration.
As in \cite{farrell2018b}, the problem on the coarsest level is solved with the SuperLU\_DIST sparse
direct solver \cite{li1999,li2003b} and uses PETSc's telescoping functionality \cite{may2016}
for improved parallel scalability.

\subsection{Verification and pressure robustness}\label{sec:svmms}
We consider the two-dimensional test case of~\cite{shih1989} to verify the implementation and to confirm that the velocity errors are independent of the Reynolds number.
The example is similar to the lid-driven cavity but with a known analytical solution.
Using either the $\PtwoPonedisc$ Scott--Vogelius or the $\PtwoPone$ Taylor--Hood element we observe the expected second order convergence of the velocity gradient and of the pressure as the mesh is refined, see Figure~\ref{fig:mms-sv-2d}.
In addition, we compare to the $\PtwoPzero$ element used in \cite{benzi2006, farrell2018b} which converges at first order only.

As motivated in the introduction, since the Taylor--Hood and the $\PtwoPzero$ element do not enforce the divergence constraint exactly (see bottom right of Figure~\ref{fig:mms-sv-2d}), the velocity error increases as the Reynolds number is increased.
This is in contrast to the solutions obtained using the Scott--Vogelius element, which are divergence-free up to solver tolerances and exhibit Reynolds-robust errors.
\begin{figure}[htbp]
  \centering
  \input{./plots/mms2d}
  \caption{Velocity and pressure error as well as $L^2$ norm of the divergence for the $\PtwoPzero$, $\PtwoPone$, and $\PtwoPonedisc$ elements for different Reynolds numbers. The $\PtwoPonedisc$ element yields velocity errors independent of the Reynolds number, while the other discretizations do not.}\label{fig:mms-sv-2d}
\end{figure}

\subsection{Two dimensional examples}
\begin{figure}[!htpb]
  \centering
  \includegraphics[height=4cm]{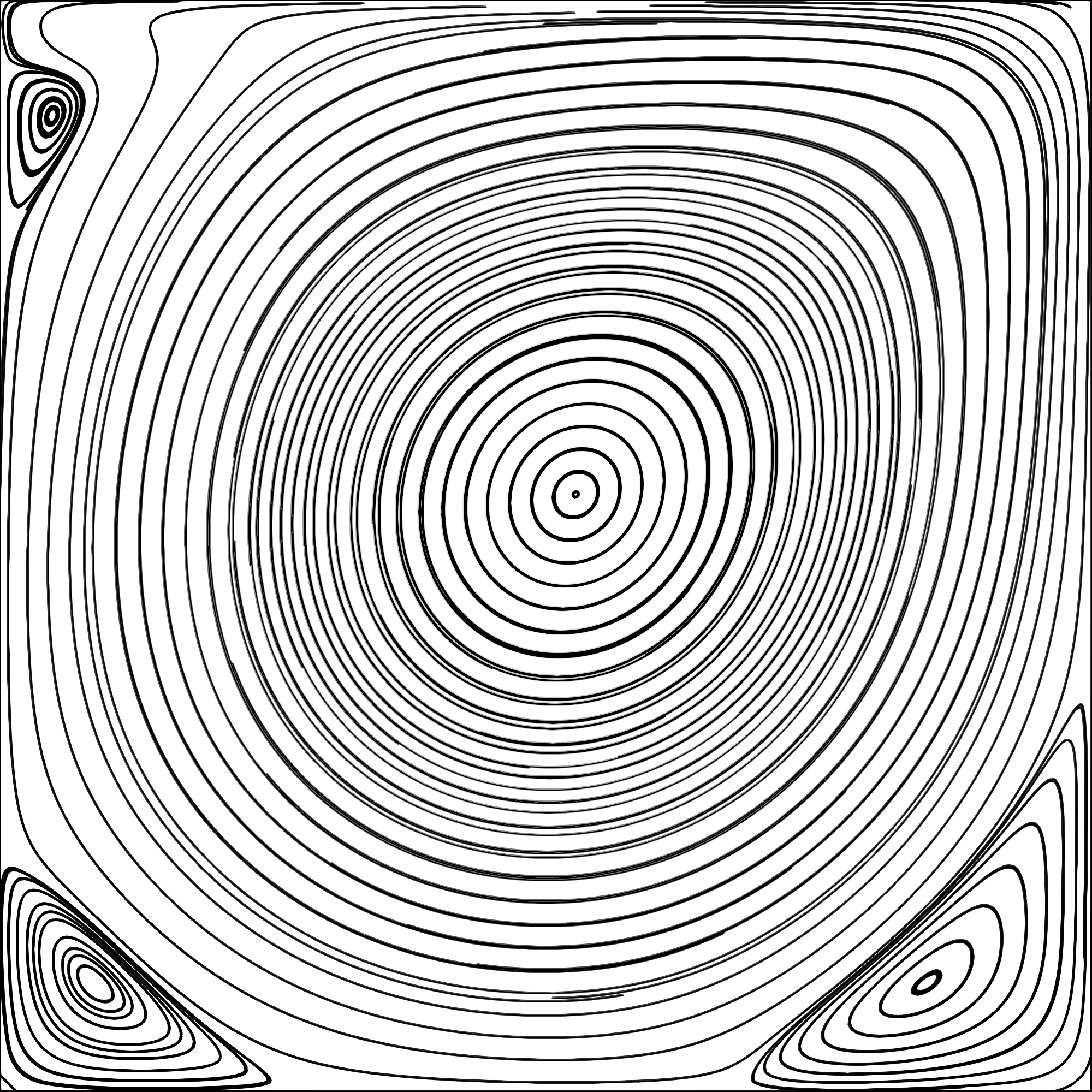}
    \caption{Streamlines for the two dimensional lid-driven cavity problem at $\Re=2\,500$. The domain is given by the $[0, 2]\times [0, 2]$ square. The boundary condition on the top is given by a horizontal velocity field $\ub(x, y)=(x^2(2-x)^2, 0)$ and the other boundaries are equipped with a no-slip condition.\label{fig:ldc2dexpl}}
\end{figure}
\begin{figure}[!htpb]
  \centering
  \includegraphics[width=0.99\textwidth]{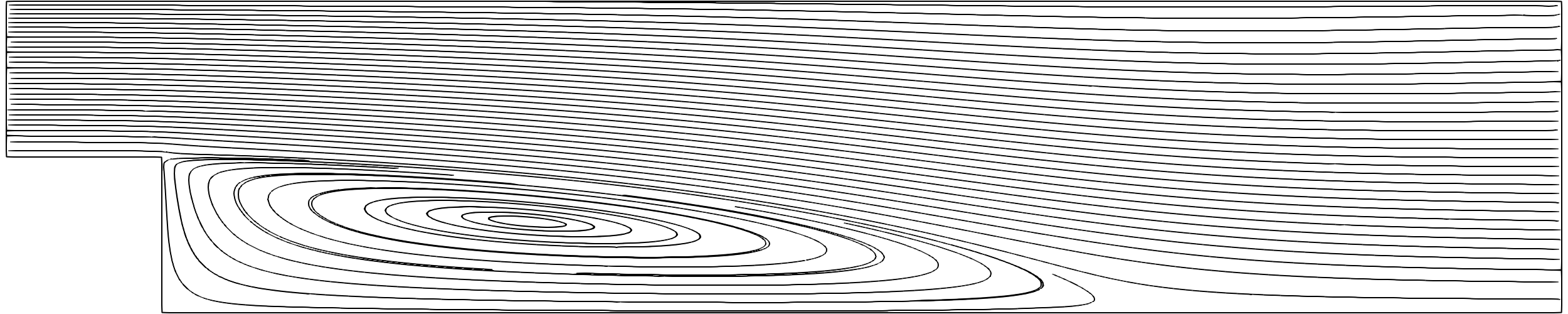}
  \caption{Streamlines for the two dimensional backwards-facing step problem at $\Re=200$. The domain is given by $([0, 10]\times [0, 2])\setminus ([0,1)\times[0,1))$. 
    The inflow condition at the top-left boundary is given by a horizontal velocity field $\ub(x,y)=(4(2-y)(y-1), 0)$, a natural outflow condition \eqref{eq:naturaloutflow} is enforced on the right and the other boundaries are equipped with a no-slip condition.\label{fig:bfs2dexpl}}
\end{figure}
We \revnew{begin by considering two classical} benchmark problems: the regularised lid-driven
cavity and backward-facing step problems, shown in Figures~\ref{fig:ldc2dexpl} and \ref{fig:bfs2dexpl}, and described in detail in 
\cite[Examples 8.1.2 and 8.1.3]{elman2014}.
We choose relative and absolute tolerances of $10^{-9}$ and $10^{-8}$ respectively for the nonlinear solver and $10^{-9}$ and $10^{-10}$ for the linear solver.
The mesh hierarchy is constructed as described in Section \ref{sec:sv-stability}.
In each case, we start with a given coarse grid, perform the specified number of uniform refinements, and then barycentrically refine each level once.

We perform continuation in the Reynolds number: we start by solving the Stokes equations, and then solve for $\Re=1$, $\Re=10$, $\Re=100$, $\Re=200$, $\Re=300$, etc., up to $\Re=10\,000$, using the previous solution as initial guess for the subsequent Newton iteration.
For the backward-facing step cases we add additional continuation steps at $\Re\in\{50, 150, 250, 350\}$.
The augmented Lagrangian parameter is set to $\gamma = 10^4$ and \revnew{we run the solver with and without the stabilization \eqref{eqn:burman-stab}.
In the former case, }the stabilization parameter is chosen as $\delta=5\cdot 10^{-3}$.

For the lid-driven cavity the coarse mesh is given by a regular $10\times 10$ triangular mesh.
Iteration counts are shown in Table~\ref{tab:svldc2d}.

\begin{table}[htbp]
\revnew{
\centering
\begin{tabular}{cc|ccccc}
\toprule
\multirow{2}{*}{Refinements} & \multirow{2}{*}{Degrees of freedom} & \multicolumn{5}{c}{Reynolds number} \\
  && 10 & 100 & 1\,000 & 5\,000 & 10\,000 \\
\midrule
\multicolumn{7}{c}{$\PtwoPonedisc$}\\
\midrule
1     & $1.70\times 10^4$ & 4.50$|$4.50 & 5.00$|$5.00 & 6.67$|$5.67 &                       18.00$|$25.00 & 25.50$|$$>$100\hphantom{.0}   \\
2     & $6.75\times 10^4$ & 4.00$|$4.00 & 4.33$|$4.33 & 5.67$|$4.33 &                       16.00$|$16.50 & 29.00$|$56.50\hphantom{$>$}\\
3     & $2.69\times 10^5$ & 4.00$|$4.00 & 4.00$|$4.00 & 4.00$|$3.00 &            10.00$|$9.50\hphantom{1} & 18.00$|$27.00\hphantom{$>$}\\
4     & $1.08\times 10^6$ & 3.00$|$3.00 & 3.33$|$3.33 & 2.67$|$2.67 & \hphantom{1}6.50$|$6.50\hphantom{1} & 11.50$|$17.50\hphantom{$>$}\\
\midrule
\multicolumn{7}{c}{$\PthreePtwodisc$}\\
\midrule
1   & $3.62\times 10^4$     & 2.50$|$2.50  & 2.67$|$2.67  & 3.33$|$2.67  & 8.00$|$8.00  & 11.50$|$15.00 \\
2   & $1.44\times 10^5$     & 2.50$|$2.50  & 2.67$|$2.67  & 2.33$|$2.00  & 5.50$|$4.50  & 9.50$|$8.50 \\
3   & $5.77\times 10^5$     & 2.00$|$2.00  & 2.67$|$2.67  & 2.00$|$1.67  & 4.00$|$3.00  & 5.50$|$5.00 \\
4   & $2.31\times 10^6$     & 2.00$|$2.00  & 2.67$|$2.67  & 2.00$|$2.00  & 2.50$|$2.00  & 5.00$|$4.00 \\
\bottomrule
\end{tabular}
\caption{Average number of outer Krylov iterations per Newton step for the
2D regularised lid-driven cavity problem using the Scott--Vogelius discretization \revnew{with and without stabilization} for $k=2$ and $k=3$.}
\label{tab:svldc2d}
}
\end{table}

Iteration counts for the backwards-facing step are shown in Table~\ref{tab:svbfs2d}.
For this case we consider a coarse mesh generated by the frontal meshing algorithm of Gmsh~\cite{geuzaine2009} consisting of 5996 triangles.

\begin{table}[htbp]
\revnew{
\centering
\begin{tabular}{cc|ccccc}
\toprule
\multirow{2}{*}{Refinements} & \multirow{2}{*}{Degrees of freedom} & \multicolumn{5}{c}{Reynolds number} \\
  && 10 & 100 & 1\,000 & 5\,000 & 10\,000 \\
\midrule
\multicolumn{7}{c}{$\PtwoPonedisc$}\\
\midrule
    1 & $4.79\times 10^5$ & 3.67$|$3.67 & 3.25$|$3.25 & 5.00$|$6.50 &            14.50$|$$>$100\hphantom{.0} & 19.50$|$$>$100            \\
    2 & $1.91\times 10^6$ & 3.67$|$3.67 & 3.25$|$3.25 & 4.00$|$4.00 & \hphantom{1}8.50$|$25.00\hphantom{$>$} & 13.50$|$$>$100            \\
    3 & $7.64\times 10^6$ & 4.50$|$4.50 & 4.33$|$4.00 & 3.00$|$3.00 & \hphantom{1}5.00$|$12.50\hphantom{$>$} & \hphantom{1}9.50$|$$>$100 \\
\midrule
\multicolumn{7}{c}{$\PthreePtwodisc$}\\
\midrule
    1 & $1.02\times 10^6$ & 2.00$|$2.00 & 2.00$|$2.00 & 2.50$|$2.50 &            4.50$|$11.00 & 6.00$|$$>$100           \\
    2 & $4.10\times 10^6$ & 2.50$|$2.50 & 2.33$|$2.33 & 1.50$|$1.50 & 2.00$|$3.50\hphantom{1} & 3.00$|$10.50            \\
    3 & $1.64\times 10^7$ & 2.50$|$2.50 & 3.33$|$3.33 & 2.00$|$2.00 & 1.50$|$2.00\hphantom{1} & 2.00$|$4.50\hphantom{1} \\
\bottomrule
\end{tabular}
\caption{Average number of outer Krylov iterations per Newton step for the
2D backwards-facing step problem using the Scott--Vogelius discretization \revnew{with and without stabilization} for $k=2$ and $k=3$.}
\label{tab:svbfs2d}
}
\end{table}

\revnew{
For both examples, we observe low and nearly flat iteration counts to $\Re=1\,000$ for both $k=2$ and $k=3$ both with and without stabilization.
It is at larger Reynolds numbers that we see differences between the four different configurations.
The combination of the $\PtwoPonedisc$ element and lack of stabilization is least robust.
Adding stabilization and increasing the polynomial degree (to additionally capture the kernel of the stabilization, cf.~Section~\ref{sec:stabilization}) results in highly robust iteration counts on the finer meshes for both the lid-driven cavity and the backwards-facing step.
}

\revnew{
\begin{remark}
    As discussed in Remark~\ref{rem:th-schur}, we expect that for the Taylor--Hood element on barycentrically-refined meshes the combination of grad-div stabilization and the presented multigrid scheme for the top-left block results in an effective preconditioner.
    To confirm this we solved the lid-driven cavity using the $[\Ptwo]^2\mathrm{-}\Pone$ and $[\Pthree]^2\mathrm{-}\Ptwo$ Taylor--Hood elements and obtained essentially identical iteration counts (up to $\pm1$) to those shown in Table~\ref{tab:svldc2d}.
\end{remark}
}

\revnew{
Clearly the patches in the multigrid relaxation considered here are significantly larger than the star patches considered in~\cite{benzi2006, farrell2018b}.
To investigate the impact this has on performance, in Table~\ref{tab:svldcruntime} we compare the runtime of the implementation in~\cite{farrell2018b} using the $\PtwoPzero$ with the $\PtwoPonedisc$ and $\PthreePtwodisc$ discretization presented here.
Since the mesh for the $\PtwoPzero$ does not require barycentric refinement, we choose a finer $16\times 16$ grid, to obtain a problem of comparable size.
As expected, the solver is more expensive for the Scott--Vogelius discretization, as it requires the larger $\MacroStar$.
However, we emphasize that the improved robustness of the $\PthreePtwodisc$ element implies that this discretization is very attractive, especially at higher Reynolds number: it converges at higher order but the computational cost is only $\sim 50\%$ bigger than that of the $\PtwoPzero$ and $\PtwoPonedisc$ discretization on the same mesh. 
}
\begin{table}[htbp]
\centering
\begin{tabular}{lc|ccccc}
\toprule
\multirow{2}{*}{Discretization} & \multirow{2}{*}{Degrees of freedom} & \multicolumn{5}{c}{Reynolds number} \\
  && 10 & 100 & 1\,000 & 5\,000 & 10\,000 \\
\midrule
$\quad \PtwoPzero$        & $6.57\times 10^5$ & 3.38  & 3.74  & 4.56  & 8.85  & 9.64\\
$\quad \PtwoPonedisc$     & $1.08\times 10^6$ & 10.39 & 10.84 & 9.06  & 18.63 & 32.31\\
$\quad \PthreePtwodisc$   & $2.31\times 10^6$ & 25.83 & 30.28 & 24.73 & 28.43 & 48.07\\
\bottomrule
\end{tabular}
\caption{Runtime per Newton step (in seconds) for the 2D regularised lid-driven cavity using the $\PtwoPzero$, $\PtwoPonedisc$, and $\PthreePtwodisc$ element pairs (with stabilization). Measured on two Intel(R) Xeon(R) Gold 5118 CPUs running 12 MPI processes.}
\label{tab:svldcruntime}
\end{table}

\revnew{
We study one final two dimensional problem.
As discussed in \cite{gauger_2019}, the nature of flows at high Reynolds number is dependent on the Helmholtz decomposition of the convection term.
Ahmed et al.~\cite{ahmed_2020} introduce the seminorm
\begin{equation}
    |f| = \sup_{0\not= v\in \left(H^1_0\right)^d, \nabla\cdot v =0} \frac{|\langle f, v\rangle|}{\|\nabla v\|_{L^2}}.
\end{equation}
This seminorm vanishes for gradient fields, i.e.~when $\langle f, v\rangle = (\nabla \phi, v)_{L^2}$ for some $\phi \in H^1_0$.
In the case when $|(u\cdot \nabla)u| = 0$ the non-linear convective term and the pressure gradient approximately balance.
These types of flows are referred to as \emph{generalized Beltrami flows} and it is particularly for these types of flows that exactly divergence-free methods were shown to outperform standard discretizations~\cite{gauger_2019}.
In particular, both of the previously studied examples are close to generalized Beltrami flows at high Reynolds numbers.
However, in general, flows may contain a significant divergence-free component in the non-linear term.
To demonstrate that the proposed preconditioner is efficient for this type of flow, we consider the Oseen problem proposed in~\cite[\S4.4]{ahmed_2020} with exact velocity $u(x, y)=(\sin(2\pi x)\sin(2\pi y), \cos(2\pi x)\cos(2\pi y))^{\top}$, pressure $p=\frac{1}4 (\cos(4\pi x)-\cos(4\pi y))$ and wind $\beta=u+(0, 1)^{\top}$ on a domain $\Omega=(0, 1)^2$.
Since this problem is linear we no longer perform continuation in the Reynolds number, and instead start the linear solver with a zero initial guess for each Reynolds number.

Iteration counts are shown in Table~\ref{tab:superposition2d} --- the results are qualitatively similar to those for the backwards-facing step: robust iteration counts both with and without stabilization to $\Re=1\,000$. Robust iteration counts to $\Re=10\,000$ are obtained when using stabilization and a cubic velocity space.
}
\begin{table}[htbp]
\revnew{
\centering
\begin{tabular}{cc|ccccc}
\toprule
\multirow{2}{*}{Refinements} & \multirow{2}{*}{Degrees of freedom} & \multicolumn{5}{c}{Reynolds number} \\
  && 10 & 100 & 1\,000 & 5\,000 & 10\,000 \\
\midrule
\multicolumn{7}{c}{$\PtwoPonedisc$}\\
\midrule
1 & $1.70\times 10^4$ & 3$|$3 & 3$|$3 & \hphantom{1}6$|$13 & 10$|$$>$100 & 10$|$$>$100 \\
2 & $6.75\times 10^4$ & 3$|$3 & 3$|$2 & \hphantom{1}6$|$11 & 15$|$$>$100 & 18$|$$>$100 \\
3 & $2.69\times 10^5$ & 3$|$3 & 2$|$2 & 4$|$7              & 15$|$$>$100 & 19$|$$>$100 \\
4 & $1.08\times 10^6$ & 3$|$3 & 2$|$2 & 3$|$5              & 12$|$$>$100 & 19$|$$>$100 \\
\midrule
\multicolumn{7}{c}{$\PthreePtwodisc$}\\
\midrule
1 & $3.62\times 10^4$ & 2$|$2 & 2$|$2 & 3$|$3 &            5$|$28 & 5$|$$>$100 \\
2 & $1.44\times 10^5$ & 2$|$2 & 2$|$2 & 3$|$3 &            4$|$21 & 6$|$$>$100 \\
3 & $5.77\times 10^5$ & 2$|$2 & 2$|$2 & 2$|$2 &            4$|$11 & 5$|$70\hphantom{$>$0} \\
4 & $2.31\times 10^6$ & 2$|$2 & 2$|$2 & 2$|$2 & 3$|$6\hphantom{1} & 4$|$25\hphantom{$>$0} \\
\bottomrule
\end{tabular}
\caption{\revnew{Number of Krylov iterations at each Reynolds number for the example in~\cite[\S4.4]{ahmed_2020} with and without stabilization for $k=2$ and $k=3$.}}
\label{tab:superposition2d}
}
\end{table}

\subsection{Three dimensional examples}
We now study three dimensional variants of the lid-driven cavity and backwards-facing step problems; these are described in detail in~\cite[\S 5.5]{farrell2018b}.
The solver tolerances are all relaxed to $10^{-8}$ in three dimensions.
We study iteration counts both with and without adding the stabilization terms in~\eqref{eqn:burman-stab}.

Results for the lid-driven cavity are shown in Table~\ref{tab:svldc3d}.
Both with and without stabilization we observe iteration counts that approximately double as the Reynolds number is increased from $\Re=10$ to $\Re=5\,000$.

We show results for the backward-facing step in Table~\ref{tab:svbfs3d}.
Without stabilization, we observe iteration counts that approximately triple over the same range of Reynolds numbers.
However, when adding stabilization, iteration counts increase significantly and blow up for very high Reynolds number.
We attribute this to the issue raised in Section~\ref{sec:stabilization}: the stabilization term itself has a large nullspace (consisting of $C^1$ vector fields) that is not captured by the relaxation induced by the $\MacroStar$ around vertices.
If we choose $k=5$ we know that a local basis for $C^1$ functions exists.
Indeed, we see in Table~\ref{tab:svbfs3d} that iteration counts for the $\PfivePfourdiscthreed$ element are significantly more robust.

\begin{table}[htbp]
\centering
\begin{tabular}{cc|ccccc}
\toprule
\multirow{2}{*}{Refinements} & \multirow{2}{*}{Degrees of freedom} & \multicolumn{5}{c}{Reynolds number} \\
  && 10 & 100 & 1\,000 & 2\,500 & 5\,000 \\
\midrule
\multicolumn{7}{c}{$k=3$ without stabilization}\\
\midrule
1   & $1.03\times 10^6$     & 3.00  & 3.67  & 3.50 & 4.00 & 5.00\\
2   & $8.22\times 10^6$     & 3.50  & 3.67  & 4.00 & 4.00 & 4.00\\
3   & $6.55\times 10^7$     & 3.00  & 3.33  & 3.50 & 3.50 & 4.00\\
\midrule
\multicolumn{7}{c}{$k=3$ with stabilization}\\
\midrule
1   & $1.03\times 10^6$     & 3.00  & 4.00  & 4.50  & 5.00 & 6.00\\
2   & $8.22\times 10^6$     & 3.50  & 4.00  & 5.50  & 6.00 & 6.50\\
3   & $6.55\times 10^7$     & 3.00  & 3.33  & 5.00  & 6.00 & 7.50\\
\bottomrule
\end{tabular}
\caption{Average number of outer Krylov iterations per Newton step for the
3D regularised lid-driven cavity problem.}
\label{tab:svldc3d}
\end{table}

\begin{table}[htbp]
\centering
\begin{tabular}{cc|ccccc}
\toprule
\multirow{2}{*}{Refinements} & \multirow{2}{*}{Degrees of freedom} & \multicolumn{5}{c}{Reynolds number} \\
  && 10 & 100 & 1\,000 & 2\,500 & 5\,000 \\
\midrule
\multicolumn{7}{c}{$k=3$ without stabilization}\\
\midrule
1	& $3.85\times 10^6$	& 4.50	& 4.33	& 5.33	& 9.00 & 15.00\\
2	& $3.06\times 10^7$	& 5.00	& 5.33	& 5.33	& 10.00& 12.00\\
\midrule
\multicolumn{7}{c}{$k=3$ with stabilization}\\
\midrule
1	& $3.85\times 10^6$	& 4.50	& 5.33	& 7.33	& 11.50 & 13.50\\
2	& $3.06\times 10^7$	& 5.00	& 6.33	& 12.50	& 14.00 & 154.00\\
  \midrule
\multicolumn{7}{c}{$k=5$ with stabilization}\\
\midrule
1   & $3.81\times 10^6$     & 1.50  & 1.67  & 2.00 &3.50 & 4.00\\
2 & $3.03 \times 10^7$ & 2.00 & 1.67 & 2.00 & 2.50 & 5.00 \\
\bottomrule
\end{tabular}
\caption{Average number of outer Krylov iterations per Newton step for the
3D backwards-facing step problem. The results for $k=5$ were obtained on a coarser mesh to have a comparable number of degrees of freedom to the case of $k=3$.}
\label{tab:svbfs3d}
\end{table}
The results in Tables~\ref{tab:svldc3d} and~\ref{tab:svbfs3d} were obtained on the ARCHER supercomputer.
To provide an impression of computational performance, the results on the twice refined mesh of Table~\ref{tab:svbfs3d} were run on 512 cores (for $k=3$, $\sim60\,000$ dofs/core) and 960 cores (for $k=5$, $\sim 31\,000$ dofs/core).
Without stabilization, each Krylov iteration takes $\sim 16$ seconds. Adding stabilization increases this to $\sim 30$ seconds (due to the larger number of non-zeros in the matrix) and increasing the order to $k=5$ increases the runtime further to $\sim 450$ seconds.
Though we would not recommend the case $k=5$ as a practical discretization due to its high computational cost and memory requirements, we include the results to demonstrate the necessity of capturing the nullspace of all singular operators in order to obtain a fully robust scheme.

\section{Summary}
The goal of this work was to develop a scalable solver for the stationary incompressible Navier--Stokes equations that exhibits both Reynolds-robust iteration counts and errors.
To achieve this goal we apply the augmented Lagrangian approach to the exactly incompressible Scott--Vogelius discretization on barycentric grids, and solve the augmented momentum block with a specialized multigrid method that exploits the barycentric structure in both relaxation and prolongation.

Dictated by inf-sup stability, the minimal polynomial degree that we require is $k=d$.
For this degree we observe robust iteration counts to $\Re \sim 1\,000$ but increase for higher Reynolds numbers.
We attribute this behaviour to the singular nature of the stabilization term employed and show that when using a sufficiently high order discretization, full robustness can be achieved.

\section*{Code availability}
The code for the Navier--Stokes solver and the numerical experiments in this paper can be found at \url{https://github.com/florianwechsung/alfi/}.
For reproducibility, all major \pkg{Firedrake} components as well as the code used to run these examples have been archived on Zenodo.
The results in Tables~\ref{tab:svldc3d}, and~\ref{tab:svbfs3d} were obtained on ARCHER, the UK national supercomputer using the code archived at~\cite{firedrake_archer, alfi_archer}. 
The results in Figure~\ref{fig:mms-sv-2d} and Tables~\ref{tab:svldc2d},~\ref{tab:svbfs2d},~\ref{tab:svldcruntime}, tand~\ref{tab:superposition2d} were obtained using the code archived at~\cite{zenodo/Firedrake-20201208.0, alfi_local}.

\ifarxiv

\else
\bibliographystyle{siamplain}
\bibliography{literature}
\fi
\end{document}